\documentclass[a4paper,pdftex,reqno,10pt]{amsart}

\allowdisplaybreaks[1]

\usepackage{geometry}

\usepackage{graphicx}
\graphicspath{{Figures/}}

\usepackage{courier}
\usepackage{times}
\usepackage{amsmath,amssymb}
\usepackage[euler-digits]{eulervm}
\usepackage{eurosym}

\usepackage{hyperref}

\theoremstyle{plain}
\newtheorem{Theorem}{Theorem}[section]

\newtheorem{Conjecture}[Theorem]{Conjecture}
\newtheorem{Corollary}[Theorem]{Corollary}
\newtheorem{Lemma}[Theorem]{Lemma}

\newenvironment{Proof}
{\begin{trivlist}\item[]{{\sc Proof.}}}{\hfill{$\square$}\noindent\end{trivlist}}

\theoremstyle{definition}
\newtheorem{Definition}[Theorem]{Definition}

\theoremstyle{remark}

\begin{document}


\title[Bounds for the minimum oriented diameter]{Bounds for the minimum oriented diameter}

\author{Sascha Kurz}

\address{Sascha Kurz\\Fakult\"at f\"ur Mathematik, Physik und Informatik\\Universit\"at Bayreuth\\Germany}
\email{sascha.kurz@uni-bayreuth.de}

\author{Martin L\"atsch}

\address{Martin L\"atsch\\Zentrum f\"ur Angewandte Informatik\\Universit\"at zu K\"oln\\Germany}
\email{laetsch@zpr.uni-koeln.de}

\begin{abstract}
  We consider the problem of finding an orientation with minimum diameter of a connected bridgeless graph. Fomin et. al.
  \cite{0981.05059} discovered a relation between the minimum oriented diameter an the size of a minimal dominating set.
  We improve their upper bound.
\end{abstract}


\maketitle

\section{Introduction}

An orientation of an undirected graph $G$ is a directed graph whose arcs correspond to assignments of directions to the edges of $G$. An Orientation $H$ of $G$ is strongly connected if every two vertices in $H$ are mutually reachable in $H$. An edge $e$ in a undirected connected graph $G$ is called a bridge if $G-e$ is not connected. A connected graph $G$ is bridgeless if $G-e$ is connected for every edge $e$, i.~e. there is no bridge in $G$.

The conditions when an undirected graph $G$ admits a strongly connected orientation are determined by Robbins in 1939 \cite{origin}. The necessary and sufficient conditions are that $G$ is connected and bridgeless. Chung et. al provided a linear-time algorithm for testing whether a graph has a strong orientation and finding one if it does \cite{linear-time}.

\begin{Definition}
  Let $\overset{\rightarrow}{G}$ be a strongly connected directed graph. By
  $diam \left(\overset{\rightarrow}{G}\right)$ we denote the diameter of $\overset{\rightarrow}{G}$.
  For a simple graph connected $G$ without bridges we define
  $$
    \overset{\longrightarrow}{diam}_{min}(G):=\min\Big\{diam\left(\overset{\rightarrow}{G}\right)\,:\,
    \overset{\rightarrow}{G}\text{ is an orientation of }G\Big\},
  $$
  which we call the minimum oriented diameter of a simple graph $G$.
  By $\gamma(G)$ we denote the smallest cardinality of a vertex cover of $G$.
\end{Definition}

We are interested in the examples $G$ which have a large minimum oriented diameter $\overset{\longrightarrow}{diam}_{min}(G)$ in dependence of its domination number $\gamma(G)$. Therefore we set
$$
  \Xi(\gamma):=\max\left\{\overset{\longrightarrow}{diam}_{min}(G)\,:\, \gamma(G)\le \gamma\text{ for $G$ being a bridgeless
  connected graph}\right\}.
$$

\noindent
The aim of this note is to prove a better upper bound on $\Xi(\gamma)$. The previously best known result \cite{0981.05059} was:

\begin{Theorem}
  \label{thm_fomin}
  $$
    \Xi(\gamma)\le 5\gamma-1.
  $$
\end{Theorem}

\noindent
Our main results are

\begin{Theorem}
  \label{thm_four_gamma}
  $$
    \Xi(\gamma)\le 4\gamma
  $$
\end{Theorem}

\noindent 
and

\begin{Conjecture}
  \label{main_conj}
  $$
    \Xi(\gamma)=\left\lceil\frac{7\gamma(G)+1}{2}\right\rceil.
  $$
\end{Conjecture}

Clearly we have $\Xi(\gamma)$ is weak monotone increasing.
At first we observe that we have $\Xi(\gamma)\ge \left\lceil\frac{7\gamma(G)+1}{2}\right\rceil$. Therefore we consider the following set of examples, where we have depicted the vertices of a possible minimal vertex cover by a filled black circle:

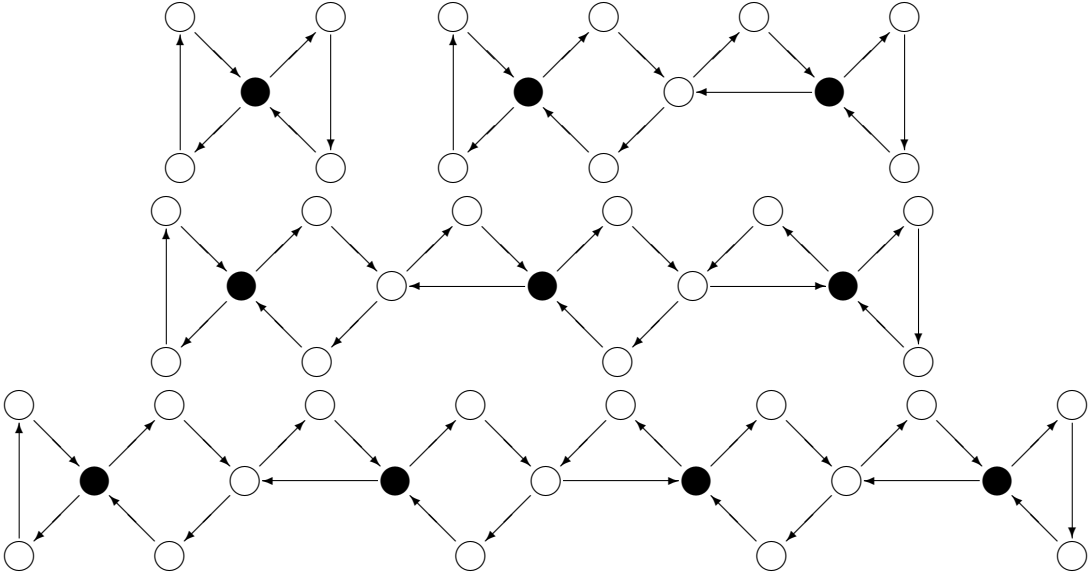
\begin{figure}[htp]
\begin{center}
  \setlength{\unitlength}{1cm}
  \begin{picture}(2.4,2.4)
    \put(0.2,0.2){\circle{0.4}}
    \put(1.2,1.2){\circle*{0.4}}
    \put(0.2,2.2){\circle{0.4}}
    \put(1,1){\vector(-1,-1){0.6}}
    \put(0.2,0.43){\vector(0,1){1.54}}
    \put(0.4,2){\vector(1,-1){0.6}}
    \put(1.4,1.4){\vector(1,1){0.6}}
    \put(2.2,2.2){\circle{0.4}}
    \put(2.2,0.2){\circle{0.4}}
    \put(2.2,1.97){\vector(0,-1){1.54}}
    \put(2.0,0.4){\vector(-1,1){0.6}}
  \end{picture}
  \quad\quad\quad
  \setlength{\unitlength}{1cm}
  \begin{picture}(6.4,2.4)
    \put(0.2,0.2){\circle{0.4}}
    \put(1.2,1.2){\circle*{0.4}}
    \put(0.2,2.2){\circle{0.4}}
    \put(1,1){\vector(-1,-1){0.6}}
    \put(0.2,0.43){\vector(0,1){1.54}}
    \put(0.4,2){\vector(1,-1){0.6}}
    \put(1.4,1.4){\vector(1,1){0.6}}
    \put(2.0,0.4){\vector(-1,1){0.6}}
    \put(2.2,2.2){\circle{0.4}}
    \put(2.2,0.2){\circle{0.4}}
    \put(2.4,2.0){\vector(1,-1){0.6}}
    \put(3,1){\vector(-1,-1){0.6}}
    \put(3.2,1.2){\circle{0.4}}
    \put(3.4,1.4){\vector(1,1){0.6}}
    \put(4.2,2.2){\circle{0.4}}
    \put(4.4,2){\vector(1,-1){0.6}}
    \put(5.2,1.2){\circle*{0.4}}
    \put(4.97,1.2){\vector(-1,0){1.54}}
    \put(5.4,1.4){\vector(1,1){0.6}}
    \put(6.2,2.2){\circle{0.4}}
    \put(6.2,0.2){\circle{0.4}}
    \put(6.2,1.97){\vector(0,-1){1.54}}
    \put(6.0,0.4){\vector(-1,1){0.6}}
  \end{picture}
  \medskip
  \setlength{\unitlength}{1cm}
  \begin{picture}(10.4,2.4)
    \put(0.2,0.2){\circle{0.4}}
    \put(1.2,1.2){\circle*{0.4}}
    \put(0.2,2.2){\circle{0.4}}
    \put(1,1){\vector(-1,-1){0.6}}
    \put(0.2,0.43){\vector(0,1){1.54}}
    \put(0.4,2){\vector(1,-1){0.6}}
    \put(1.4,1.4){\vector(1,1){0.6}}
    \put(2.0,0.4){\vector(-1,1){0.6}}
    \put(2.2,2.2){\circle{0.4}}
    \put(2.2,0.2){\circle{0.4}}
    \put(2.4,2.0){\vector(1,-1){0.6}}
    \put(3,1){\vector(-1,-1){0.6}}
    \put(3.2,1.2){\circle{0.4}}
    \put(3.4,1.4){\vector(1,1){0.6}}
    \put(4.2,2.2){\circle{0.4}}
    \put(4.4,2){\vector(1,-1){0.6}}
    \put(5.2,1.2){\circle*{0.4}}
    \put(4.97,1.2){\vector(-1,0){1.54}}
    \put(5.4,1.4){\vector(1,1){0.6}}
    \put(6.0,0.4){\vector(-1,1){0.6}}
    \put(6.2,2.2){\circle{0.4}}
    \put(6.2,0.2){\circle{0.4}}
    \put(6.4,2.0){\vector(1,-1){0.6}}
    \put(7,1){\vector(-1,-1){0.6}}
    \put(7.2,1.2){\circle{0.4}}
    \put(8,2){\vector(-1,-1){0.6}}
    \put(8.2,2.2){\circle{0.4}}
    \put(9,1.4){\vector(-1,1){0.6}}
    \put(9.2,1.2){\circle*{0.4}}
    \put(7.43,1.2){\vector(1,0){1.54}}
    \put(9.4,1.4){\vector(1,1){0.6}}
    \put(10.2,2.2){\circle{0.4}}
    \put(10.2,0.2){\circle{0.4}}
    \put(10.2,1.97){\vector(0,-1){1.54}}
    \put(10.0,0.4){\vector(-1,1){0.6}}
  \end{picture}
  \medskip
  \setlength{\unitlength}{1cm}
  \begin{picture}(14.4,2.4)
    \put(0.2,0.2){\circle{0.4}}
    \put(1.2,1.2){\circle*{0.4}}
    \put(0.2,2.2){\circle{0.4}}
    \put(1,1){\vector(-1,-1){0.6}}
    \put(0.2,0.43){\vector(0,1){1.54}}
    \put(0.4,2){\vector(1,-1){0.6}}
    \put(1.4,1.4){\vector(1,1){0.6}}
    \put(2.0,0.4){\vector(-1,1){0.6}}
    \put(2.2,2.2){\circle{0.4}}
    \put(2.2,0.2){\circle{0.4}}
    \put(2.4,2.0){\vector(1,-1){0.6}}
    \put(3,1){\vector(-1,-1){0.6}}
    \put(3.2,1.2){\circle{0.4}}
    \put(3.4,1.4){\vector(1,1){0.6}}
    \put(4.2,2.2){\circle{0.4}}
    \put(4.4,2){\vector(1,-1){0.6}}
    \put(5.2,1.2){\circle*{0.4}}
    \put(4.97,1.2){\vector(-1,0){1.54}}
    \put(5.4,1.4){\vector(1,1){0.6}}
    \put(6.0,0.4){\vector(-1,1){0.6}}
    \put(6.2,2.2){\circle{0.4}}
    \put(6.2,0.2){\circle{0.4}}
    \put(6.4,2.0){\vector(1,-1){0.6}}
    \put(7,1){\vector(-1,-1){0.6}}
    \put(7.2,1.2){\circle{0.4}}
    \put(8,2){\vector(-1,-1){0.6}}
    \put(8.2,2.2){\circle{0.4}}
    \put(9,1.4){\vector(-1,1){0.6}}
    \put(9.2,1.2){\circle*{0.4}}
    \put(7.43,1.2){\vector(1,0){1.54}}
    \put(9.4,1.4){\vector(1,1){0.6}}
    \put(10.0,0.4){\vector(-1,1){0.6}}
    \put(10.2,2.2){\circle{0.4}}
    \put(10.2,0.2){\circle{0.4}}
    \put(10.4,2.0){\vector(1,-1){0.6}}
    \put(11,1){\vector(-1,-1){0.6}}
    \put(11.2,1.2){\circle{0.4}}
    \put(11.4,1.4){\vector(1,1){0.6}}
    \put(12.2,2.2){\circle{0.4}}
    \put(12.4,2){\vector(1,-1){0.6}}
    \put(13.2,1.2){\circle*{0.4}}
    \put(12.97,1.2){\vector(-1,0){1.54}}
    \put(13.4,1.4){\vector(1,1){0.6}}
    \put(14.2,2.2){\circle{0.4}}
    \put(14.2,0.2){\circle{0.4}}
    \put(14.2,1.97){\vector(0,-1){1.54}}
    \put(14.0,0.4){\vector(-1,1){0.6}}
  \end{picture}
  \caption{Examples with large minimum oriented diameter in dependence of the domination number $\gamma(G)$.}
  \label{fig_the_bad_examples}
\end{center}
\end{figure}

If we formalize this construction of graphs $G$, which is depicted for $\gamma(G)=\gamma=1,2,3,4$ we obtain examples which attain the proposed upper bound $\left\lceil\frac{7\gamma(G)+1}{2}\right\rceil$ for all $\gamma\in\mathbb{N}$. In the following we always depict vertices in a given vertex cover by a filled circle.

\subsection{Related results}
Instead of an upper bound of $\overset{\longrightarrow}{diam}_{min}(G)$ in dependence of $\gamma(G)$ on
is also interested in an upper bound in dependence of the diameter $diam(G)$. Here the best known result is given by \cite{0311.05115}:

\begin{Theorem}(Chv\'atal and Thomassen, 1978)
 Let $f(d)$ denote the best upper bound on $\overset{\longrightarrow}{diam}_{min}(G)$ where $d=diam(G)$ and $G$ is connected and bridgeless.\\
  If $G$ is a connected bridgeless graph then we have
  $$
    \frac{1}{2}diam(G)^2 +diam(G)\le f(d) \le 2\cdot diam(G)
    \cdot(diam(G)+1).
  $$
\end{Theorem}

In \cite{0311.05115} it was also shown that we have $f(2)=6$. Examples achieving this upper bound are given by the Petersen graph and by the graph obtained from $K_4$ by subdividing the three edges incident to one vertex. Recently in \cite{diameter_3} $9\le f(3)\le 11$ was shown.

The oriented diameter is trivially at least the diameter. Graphs where equality holds are said to be \textit{tight}. In \cite{0941.05035} some Cartesian products of graphs are shown to be tight. For $n\ge 4$ the $n$-cubes are tight \cite{0704.05030}. The discrete tori $C_n\times C_m$ which are tight are completely determined in \cite{0990.05041}.

The origin of this problem goes back to 1938, where Robbins \cite{origin} proves that a graph $G$ has a strongly connected orientation if and only if $G$ has no cut-edge. As an application one might think of  making streets of a city one-way or building a communication network with links that are reliable only in one direction.

There is a huge literature on the minimum oriented diameter for special graph classes, see i.e. \cite{1064.05071,0846.05038,0861.05029,0886.05060,0945.05018,0946.05036,0998.05016,1104.05024,0613.05024}.

From the algorithmic point of view the following result is known \cite{0311.05115}:

\begin{Theorem}
  The problem whether $\overset{\longrightarrow}{diam}_{min}(G)\le 2$ is $\mathcal{NP}$-hard for a
  given graph $G$.
\end{Theorem}

We remark that the proof is based on a transformation to the problem whether a hypergraph of rank $3$ is two-colorable.

\section{Preliminaries}
A vertex set $D\subseteq V(G)$ of a graph $G$ is said to be a dominating set of $G$ if for every vertex $u\in V(G)\backslash D$ there is a vertex $v\in D$ such that $\{u,v\}\in E(G)$. The minimum cardinality of a dominating set of a graph $G$ is denoted by $\gamma(G)$. If $P$ is a path we denote by $|P|$ its length which equals the number of its edges. An elementary cycle $C$ of a graph $G=(V,E)$ is a list $[v_0,\dots,v_k]$ of vertices in $V$, where $v_0=v_k$, $|\{v_0,\dots,v_k-1\}|=k$ and $\{v_i,v_{i+1}\}\in E$ for $0\le i<k$. Similarly $|C|$ denotes the length of $C$ which equals the number of its edges and vertices. For other not explicitly mention graph-theoretic terminology we refer the reader to \cite{diestel} for the basic definitions.

Our strategy to prove bounds on $\Xi(\gamma)$ is to apply some transformations on bridgeless connected graphs attaining $\Xi(\gamma)$ to obtain some structural results. Instead of considering graphs $G$ from now on we will always consider pairs $(G,D)$, where $D$ is a dominating set of $G$.

\begin{Definition}
  For a graph $G$ and a dominating set $D$ of $G$ we call $\{u,v\}\subseteq V(G)\backslash D$ an isolated triangle if
  there exists an $w\in D$ such that all neighbors of $u$ and $v$ are contained in $\{u,v,w\}$ and $\{u,v\}\in E(G)$.
  We say that the isolated triangle is associated with $w\in D$.
\end{Definition}

\begin{Definition}\label{Def_first_standard}
  A pair $(G,D)$ is in first standard form if
  \begin{itemize}
    \item[(1)] $G=(V,E)$ is a connected simple graph without a bridge,
    \item[(2)] $D$ is a dominating set of $G$ with $|D|=\gamma(G)$,
    \item[(3)] for $u,v\in D$ we have $\{u,v\}\notin E$,
    \item[(4)] for each $u\in V\backslash D$ there exists exactly one $v\in D$ with $\{u,v\}\in E$, and
    \item[(5)] $G$ is edge-minimal, meaning one can not delete an edge in $G$ without creating a bridge,
               destroying the connectivity or destroying the property of $D$ being a dominating set,
    \item[(6)] for $|D|=\gamma(G)\ge 2$ every vertex in $D$ is associated with exactly one isolated triangle
               and for $|D|=\gamma(G)=1$ the vertex in $D$ is associated with exactly two isolated triangles.
  \end{itemize}
\end{Definition}

\begin{Lemma}
  \label{lemma_first_standard}
  $$
    \Xi(\gamma)=\max\left\{\overset{\longrightarrow}{diam}_{min}(G)\,:\, |D|\le\gamma,\, (G,D)\text{ is in
    first standard form}\right\}.
  $$
\end{Lemma}
\begin{Proof}
  For a given $\gamma\in\mathbb{N}$ we start with a bridgeless connected graph $G'$ attaining
  $\Xi(\gamma)=\overset{\longrightarrow}{diam}_{min}(G')$ and minimum domination number $\gamma(G')$.
  Let $D'$ be an arbitrary dominating set of $G'$ fulfilling $|D'|=\gamma(G')$. Our aim is to apply
  some graph transformations onto $(G',D')$ to obtain a pair $(G,D)$ in first standard form fulfilling
  $\overset{\longrightarrow}{diam}_{min}(G)\ge \overset{\longrightarrow}{diam}_{min}(G')$ and $|D|\le |D'|$.

  At the start conditions (1) and (2) are fulfilled. If there is an edge $e$
  between two nodes of $D$ then we recursively apply the following graph transformation until there
  exists no such edge:
  \begin{center}
    \setlength{\unitlength}{1cm}
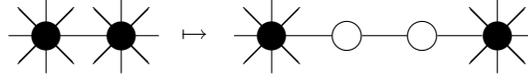
    \begin{picture}(6.4,2)
      \put(0.2,1){\circle*{0.4}}
      \put(1.2,1){\circle*{0.4}}
      \put(0.2,1){\line(1,0){1}}
      \put(0.2,1){\line(0,1){0.5}}
      \put(0.2,1){\line(-1,1){0.37}}
      \put(0.2,1){\line(-1,0){0.5}}
      \put(0.2,1){\line(-1,-1){0.37}}
      \put(0.2,1){\line(0,-1){0.5}}
      \put(0.2,1){\line(1,1){0.37}}
      \put(0.2,1){\line(1,-1){0.37}}
      \put(1.2,1){\line(0,1){0.5}}
      \put(1.2,1){\line(-1,1){0.37}}
      \put(1.2,1){\line(1,0){0.5}}
      \put(1.2,1){\line(-1,-1){0.37}}
      \put(1.2,1){\line(0,-1){0.5}}
      \put(1.2,1){\line(1,1){0.37}}
      \put(1.2,1){\line(1,-1){0.37}}
      \put(2.0,0.93){$\mapsto$}
      \put(3.2,1){\circle*{0.4}}
      \put(4.2,1){\circle{0.4}}
      \put(5.2,1){\circle{0.4}}
      \put(6.2,1){\circle*{0.4}}
      \put(3.2,1){\line(1,0){0.8}}
      \put(4.4,1){\line(1,0){0.6}}
      \put(5.4,1){\line(1,0){0.6}}
      \put(3.2,1){\line(0,1){0.5}}
      \put(3.2,1){\line(-1,1){0.37}}
      \put(3.2,1){\line(-1,0){0.5}}
      \put(3.2,1){\line(-1,-1){0.37}}
      \put(3.2,1){\line(0,-1){0.5}}
      \put(3.2,1){\line(1,1){0.37}}
      \put(3.2,1){\line(1,-1){0.37}}
      \put(6.2,1){\line(0,1){0.5}}
      \put(6.2,1){\line(-1,1){0.37}}
      \put(6.2,1){\line(1,0){0.5}}
      \put(6.2,1){\line(-1,-1){0.37}}
      \put(6.2,1){\line(0,-1){0.5}}
      \put(6.2,1){\line(1,1){0.37}}
      \put(6.2,1){\line(1,-1){0.37}}
    \end{picture}
  \end{center}

  If there exists a node $v\in V\backslash D$ with at least $r\ge 2$ neighbors $d_1,\dots,d_r$ in
  $D$ then we replace the edge $(v,d_i) \ \ i=2,\ldots,r$ with a path of length $2$.
  We iterate this until case (4) is fulfilled. In Figure \ref{Fig_punkt_4} we have depicted the graph
  transformation for $r=2,3$.

  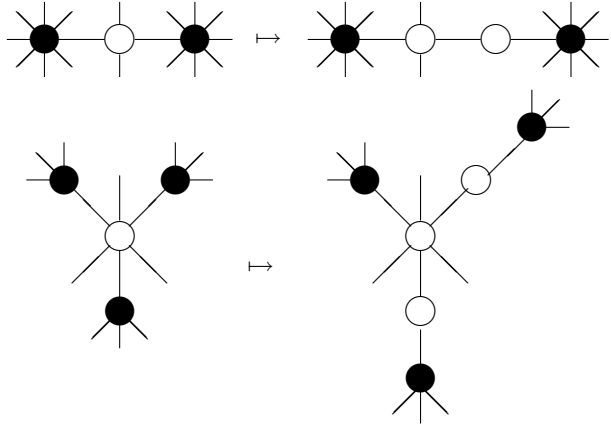
\begin{figure}[htb]
    \begin{center}
      \setlength{\unitlength}{1.0cm}
      \begin{picture}(7.4,2)
        \put(0.2,1){\circle*{0.4}}
        \put(1.2,1){\circle{0.4}}
        \put(2.2,1){\circle*{0.4}}
        \put(0.2,1){\line(1,0){0.8}}
        \put(1.4,1){\line(1,0){0.8}}
        \put(0.2,1){\line(0,1){0.5}}
        \put(0.2,1){\line(-1,1){0.37}}
        \put(0.2,1){\line(-1,0){0.5}}
        \put(0.2,1){\line(-1,-1){0.37}}
        \put(0.2,1){\line(0,-1){0.5}}
        \put(0.2,1){\line(1,1){0.37}}
        \put(0.2,1){\line(1,-1){0.37}}
        \put(1.2,1.2){\line(0,1){0.3}}
        \put(1.2,0.8){\line(0,-1){0.3}}
        \put(2.2,1){\line(0,1){0.5}}
        \put(2.2,1){\line(-1,1){0.37}}
        \put(2.2,1){\line(1,0){0.5}}
        \put(2.2,1){\line(-1,-1){0.37}}
        \put(2.2,1){\line(0,-1){0.5}}
        \put(2.2,1){\line(1,1){0.37}}
        \put(2.2,1){\line(1,-1){0.37}}
        \put(3.0,0.93){$\mapsto$}
        \put(4.2,1){\circle*{0.4}}
        \put(5.2,1){\circle{0.4}}
        \put(6.2,1){\circle{0.4}}
        \put(7.2,1){\circle*{0.4}}
        \put(4.2,1){\line(1,0){0.8}}
        \put(5.4,1){\line(1,0){0.6}}
        \put(6.4,1){\line(1,0){0.6}}
        \put(4.2,1){\line(0,1){0.5}}
        \put(4.2,1){\line(-1,1){0.37}}
        \put(4.2,1){\line(-1,0){0.5}}
        \put(4.2,1){\line(-1,-1){0.37}}
        \put(4.2,1){\line(0,-1){0.5}}
        \put(4.2,1){\line(1,1){0.37}}
        \put(4.2,1){\line(1,-1){0.37}}
        \put(5.2,1.2){\line(0,1){0.3}}
        \put(5.2,0.8){\line(0,-1){0.3}}
        \put(7.2,1){\line(0,1){0.5}}
        \put(7.2,1){\line(-1,1){0.37}}
        \put(7.2,1){\line(1,0){0.5}}
        \put(7.2,1){\line(-1,-1){0.37}}
        \put(7.2,1){\line(0,-1){0.5}}
        \put(7.2,1){\line(1,1){0.37}}
        \put(7.2,1){\line(1,-1){0.37}}
      \end{picture}

      \smallskip

      \setlength{\unitlength}{1.0cm}
      \begin{picture}(9,4)
        \put(2,2.5){\circle{0.4}}
        \put(1.87,2.63){\line(-1,1){0.61}}
        \put(2.13,2.63){\line(1,1){0.61}}
        \put(2,2.3){\line(0,-1){0.8}}
        \put(2,1.5){\circle*{0.4}}
        \put(2.74,3.24){\circle*{0.4}}
        \put(1.26,3.24){\circle*{0.4}}
        \put(2,1.5){\line(0,-1){0.5}}
        \put(2,1.5){\line(-1,-1){0.37}}
        \put(2,1.5){\line(1,-1){0.37}}
        \put(1.26,3.24){\line(-1,0){0.5}}
        \put(1.26,3.24){\line(0,1){0.5}}
        \put(1.26,3.24){\line(-1,1){0.37}}
        \put(2.74,3.24){\line(1,0){0.5}}
        \put(2.74,3.24){\line(0,1){0.5}}
        \put(2.74,3.24){\line(1,1){0.37}}
        \put(2,2.7){\line(0,1){0.6}}
        \put(2.13,2.37){\line(1,-1){0.5}}
        \put(1.87,2.37){\line(-1,-1){0.5}}
        \put(3.7,2){$\mapsto$}
        \put(6,2.5){\circle{0.4}}
        \put(5.87,2.63){\line(-1,1){0.61}}
        \put(6.13,2.63){\line(1,1){0.5}}
        \put(6,2.3){\line(0,-1){0.61}}
        \put(6,1.5){\circle{0.4}}
        \put(6,0.6){\circle*{0.4}}
        \put(6.74,3.24){\circle{0.4}}
        \put(5.26,3.24){\circle*{0.4}}
        \put(7.48,3.94){\circle*{0.4}}
        \put(6,1.25){\line(0,-1){0.61}}
        \put(6,0.5){\line(0,-1){0.5}}
        \put(6,0.5){\line(-1,-1){0.37}}
        \put(6,0.5){\line(1,-1){0.37}}
        \put(5.26,3.24){\line(-1,0){0.5}}
        \put(5.26,3.24){\line(0,1){0.5}}
        \put(5.26,3.24){\line(-1,1){0.37}}
        \put(7.48,3.94){\line(1,0){0.5}}
        \put(7.48,3.94){\line(0,1){0.5}}
        \put(7.48,3.94){\line(1,1){0.37}}
        \put(6.9,3.3){\line(1,1){0.61}}
        \put(6,2.7){\line(0,1){0.6}}
        \put(6.13,2.37){\line(1,-1){0.5}}
        \put(5.87,2.37){\line(-1,-1){0.5}}
      \end{picture}
    \end{center}
    \caption{Graph transformation to fulfill condition $(4)$ of Definition
    \ref{Def_first_standard}}\label{Fig_punkt_4}
  \end{figure}

  \bigskip

  So after a finite number of transformation we have constructed a pair $(G,D)$ which fulfills conditions
  (1), (3), (4) of the first standard form where $D$ is a dominating set of $G$ and $(G,D)$ also fulfills
  $$
    \gamma(G)\le |D|\le |D'|=\gamma(G')
  $$
  and
  $$
   \infty>\overset{\longrightarrow}{diam}_{min}(G)\ge\overset{\longrightarrow}{diam}_{min}(G').
  $$

  To additionally fulfill condition (5) of the first standard form we only need to delete the
  controversial edges. If $\gamma(G)<|D|\le\gamma(G')$ we would have a contradiction to the minimality of
  $\gamma(G')$. Since adding isolated triangles to does not contradict with the other properties and also does not
  decrease the minimum oriented property we can assume that every vertex of $D$ is associated with enough isolated
  triangles. For two vertices $x$ and $y$ in two different isolated triangles being associated with the same vertex
  $w\in D$ we have $d(x,y)\le 4$ in every strongly connected orientation. Thus we can delete some isolated triangles
  to achieve the stated number of isolated triangles for every vertex in the dominating set $D$. Finally we have a pair
  $(G,D)$ in first standard form.
\end{Proof}

So in order to prove bounds on $\Xi(\gamma)$ we can restrict ourselves on pairs $(G,D)$ in first standard form. Due to Theorem \ref{thm_fomin} we can assume $\gamma(G)=|D|\ge 2$ both for the proof of Theorem \ref{thm_four_gamma} and also for Conjecture \ref{main_conj}.

\begin{Corollary}\label{cor_function_f}
  If $(G,D)$ is a pair in first standard form then we have
  \begin{itemize}
    \item[(i)] for all $u,v\in D$ the distance fulfills $d(u,v)\ge 3$ and
    \item[(ii)] for all $u\in V(G)\backslash D$ there exists exactly one $f(u)\in D$ with $\{u,f(u)\}\in E(G)$.
  \end{itemize}
\end{Corollary}

Let $G$ be a bridgeless connected undirected graph, $D$ be a dominating set of $G$ and $H$ be a strongly connected orientation of $G$. By $diam_i(H,D)$ we denote 
$$
  \max\left\{d_H(u,v)\,:\,\Big|\{u,v\}\cap (V(H)\backslash D)\Big|=i\right\}.
$$
Clearly we have $diam(H)=\max\Big\{diam_0(H,D),diam_1(H,D),diam_2(H,D)\Big\}$. Now we refine a lemma from \cite{0981.05059}:

\begin{Lemma}
  \label{lemma_subgraph}
  Let $G'$ and $G$ be bridgeless connected graphs such that $G$ is a subgraph of $G'$ and $D$ is a dominating set of both $G'$ and
  $G$. Then for every strongly connected orientation $H$ of $G$ there is an orientation $H'$ of $G'$ such that
  $$
    diam(H')\le \max\Big\{diam_0(H,D)+4,diam_1(H,D)+2,diam_2(H,D)\Big\}.
  $$ 
\end{Lemma}
\begin{Proof}
  (We rephrase most of the proof from \cite{0981.05059}.) We adopt the direction of the edges from $H$ to $H'$. For the remaining
  edges we consider connected components $Q$ of $G'\backslash V(G)$ and direct some edges having ends in $Q$ as follows.

  If $Q$ consists of one vertex $x$ then $x$ is adjacent to at least one vertex $u$ in $D$ and to another vertex $v\neq u$ (the
  graph $G$ is bridgeless and $D$ is a dominating set). If also $v$ is an element of $D$ then we direct one edge from $x$ and
  the second edge towards $x$. Otherwise $v$ is in $V\backslash D$. In this case we direct the edges $[x,u]$ and $[v,x]$ in the
  same direction as the edge $[f(v),v]$. If there are more edges incident with $x$ (in both cases) we direct them arbitrarily.
  Then, we have assured the existence of vertices $u',v'\in D$ such that $d_{H'}(x,v')\le 2$ and $d_{H'}(u',x)\le 2$.

  Suppose that there are at least two vertices in the connected component $Q$. Choose a spanning tree $T$ in this component
  rooted in a vertex $v$. We orient edges of this tree as follows: If a vertex $x$ of the tree has odd distance from $v$, then
  we orient all the tree edges adjacent to $x$ from $x$ outwards. Also, for every such vertex $x$ we orient the edges between $x$
  and $V(G)$ towards $x$ if the distance from $v$ on the tree is even, and towards $V(G)$ otherwise, see Figure 1 in
  \cite{0981.05059}. The rest of the edges in the connected component $Q$ are oriented arbitrarily.

  In such an orientation $H'$, for every vertex $x\in Q$ there are vertices $u,v\in D$ such that $d_{H'}(x,v)\le 2$ and
  $d_{H'}(u,x)\le 2$. Therefore, for every $x,y\in V(G')$ the distance between $x$ and $y$ in $H'$ is at most
  $$
    \max\Big\{diam_0(H,D)+4,diam_1(H,D)+2,diam_2(H,D)\Big\}.
  $$
\end{Proof}

Due to the isolated triangles being associated with the vertices of the dominating set $D$, for every pair $(G,D)$ in
first standard form, there exists an orientation $H$ of $G$ such that
\begin{equation}
  \label{eq_minimal_orientation}
  \overset{\longrightarrow}{diam}_{min}(G)=diam(H)=\max\Big\{diam_0(H,D)+4,diam_1(H,D)+2,diam_2(H,D)\Big\}.
\end{equation}
If we say that $H$ is an optimal or an minimal orientation of $(G,D)$ we mean an orientation that fulfills Equation \ref{eq_minimal_orientation}.

In \cite{0981.05059} the authors have described a nice construction to obtain such a subgraph $G$ for a given bridgeless connected graph $G'$ fulfilling $|V(G)|\le 5\cdot\gamma(G')-4$:

For $\gamma(G')=1$ we may simply choose the single vertex in $D$ as our subgraph $D$. Now we assume $|D|=\gamma(G')\ge 2$. Iteratively, we construct a tree $T_k$ for $k=1,\dots,|D|$. The tree $T_1$ is composed by one vertex $x_1$ in $D$. To construct $T_{k+1}$ from $T_k$ we find a vertex $x_{k+1}$ in $D\backslash V(T_k)$ with minimum distance to $T_k$. The tree $T_{k+1}$ is the union of $T_k$ with a shortest path from $x_{k+1}$ to $T_k$. Since $D$ is a dominating set this path has length at most $3$. We say that the edges of this path are \textit{associated} with $x_{k+1}$. At the last step we obtain a dominating tree $T$ with $D\subseteq T$ and with $|V(T)|\le 2(|D|-1)+|D|$.

In order to transform $T$ in a bridgeless connected graph we construct a sequence of subgraphs $G_k$ for $k=1,\dots,|D|$. We say that $x_j\in D$ is \textit{fixed} in $G_k$ if no edge associated with $x_j$ is a bridge in $G_k$. Notice that $x_1$ is fixed in $T$ because it does not have any associated edge.

We set $G_1=T$. Assume we have constructed the subgraph $G_k$. If $x_{k+1}$ is already fixed in $G_k$ we set $G_{k+1}=G_k$. If $x_{k+1}$ is not fixed in $G_k$ we add a subgraph $M$ to $G_k$ to obtain $G_{k+1}$.

Let $P_k$ be the path added to $T_k$ to obtain $T_{k+1}$. We only consider the case where $P_k$ has length three. The other cases cane be done similarly. Let us assume that $P_k$ is given by $P_k=(x_{k+1},u,v,x_j)$ with $u,v\notin D$, and $x_j\in D$, $j\le k$. Moreover let us denote the edges of $P_k$ by $e$, $e'$ and $e''$. If we remove all edges $e$, $e'$, $e''$ of $P_k$ from $T$ we obtain four subtrees $T^1$, $T^2$, $T^3$ and $T^4$ containing $x_{k+1}$, $u$, $v$ and $x_j$, respectively.

Among all shortest path in $G'\backslash e$ connecting $T^1$ with $T^2\cup T^3\cup T^4$ we select $P$ as one whose last vertex belongs to $T^i$ with $i$ maximum. Among all shortest path in $G'\backslash e''$ connecting $T^4$ with $T^1\cup T^2\cup T^3$ we select $Q$ as one whose first vertex belongs to $T^i$ with $i$ minimum. Let $R$ be any shortest path in $G'\backslash e'$ connecting $T^3\cup T^4$ with $T^1\cup T^2$.

Since $G'$ is a bridgeless connected graph the paths $P$, $Q$, $R$ exist. Since $D\subseteq V(T)$ and the set $D$ is a dominating set, the length of paths $P$, $Q$ and $R$ is at most $3$. Moreover, if the length of $P$ is three its end vertices belong to $D$. The same holds for the paths $Q$ and $R$.

The definition of $M$ is given according to the following cases. If the last vertex of $P$ belongs to $T^4$ we define $M=P$. If the last vertex of $P$ belongs to $T^3$ or it belongs to $T^2$ and the first vertex of $Q$ belongs to $T^2$ we define $M=P\cup Q$. If none of the previous cases hold the first vertex of $R$ belongs to $T^2$ and the last one belongs to $T^3$. We define $M=P\cup Q\cup R$.

For the analysis that $|V(G_{|D|})|\le 5\cdot\gamma(G')-4$ we refer to \cite{0981.05059}.

\bigskip

Since a shortest path does contain every vertex at most once, we can combine the above described construction of a subgraph with Lemma \ref{lemma_subgraph} to obtain the bound $\Xi(\gamma)\le 5\gamma -1$.

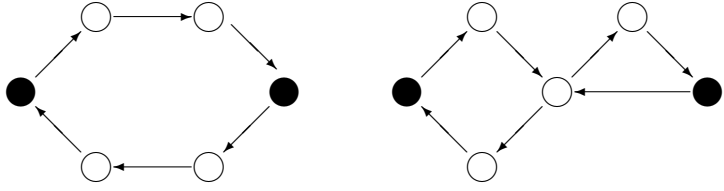
\begin{figure}[htp]
  \begin{center}
    \setlength{\unitlength}{1cm}
    \begin{picture}(3.9,2.4)
      \put(0.2,1.2){\circle*{0.4}}
      \put(0.4,1.4){\vector(1,1){0.6}}
      \put(1.2,2.2){\circle{0.4}}
      \put(1.2,0.2){\circle{0.4}}
      \put(1.0,0.4){\vector(-1,1){0.6}}
      \put(1.43,2.2){\vector(1,0){1.04}}
      \put(2.47,0.2){\vector(-1,0){1.04}}
      \put(2.7,2.2){\circle{0.4}}
      \put(2.7,0.2){\circle{0.4}}
      \put(3.5,1.0){\vector(-1,-1){0.6}}
      \put(3.0,2.1){\vector(1,-1){0.6}}
      \put(3.7,1.2){\circle*{0.4}}
    \end{picture}
    \quad\quad\quad
    \setlength{\unitlength}{1cm}
    \begin{picture}(4.4,2.4)
      \put(0.2,1.2){\circle*{0.4}}
      \put(0.4,1.4){\vector(1,1){0.6}}
      \put(1.0,0.4){\vector(-1,1){0.6}}
      \put(1.2,2.2){\circle{0.4}}
      \put(1.2,0.2){\circle{0.4}}
      \put(1.4,2.0){\vector(1,-1){0.6}}
      \put(2,1){\vector(-1,-1){0.6}}
      \put(2.2,1.2){\circle{0.4}}
      \put(2.4,1.4){\vector(1,1){0.6}}
      \put(3.2,2.2){\circle{0.4}}
      \put(3.4,2){\vector(1,-1){0.6}}
      \put(4.2,1.2){\circle*{0.4}}
      \put(3.97,1.2){\vector(-1,0){1.54}}
    \end{picture}
    \caption{The two possible subgraphs for $\gamma(G)=2$.}
    \label{fig_gamma_2}
  \end{center}
\end{figure}

\begin{Lemma}
  \label{lemma_small_exact_values}
  $$
    \Xi(1)=4\text{ and }\Xi(2)=8.
  $$
\end{Lemma}
\begin{Proof}
  At first we observe that the examples from Figure \ref{fig_the_bad_examples} give $\Xi(1)\ge 4$ and $\Xi(2)\ge 8$. For
  the other direction let $(G,D)$ be a pair in first standard form attaining 
  $\overset{\longrightarrow}{diam}_{min}(G)=\Xi(\gamma(G))$. For $\gamma=\gamma(G)=1$ we have $|D|=1$, choose the single
  vertex of $D$ as a subgraph and apply Lemma \ref{lemma_subgraph}. Going through the cases of the above described subgraph
  construction for $\gamma=\gamma(G)=2$ we obtain up to symmetry the two possibilities given in Figure \ref{fig_gamma_2}.
  By $H$ be denote the depicted corresponding orientation of the edges. Since in both cases we have $diam_0(H,D)\le 4$
  and $diam_1(H,D),diam_2(H,D)\le 5$ we can apply Lemma \ref{lemma_subgraph} to obtain the stated result.
\end{Proof}

With Lemma \ref{lemma_subgraph} in mind we would like to restrict our investigations on bridgeless connected subgraphs containing the dominating set.

\begin{Definition}
  \label{def_minimal_subgraph}
  For a pair $(G',D)$ in first standard form we call $G$ a minimal subgraph of $(G',D)$, if
  \begin{itemize}
    \item[(1)] $G$ is a subgraph of $G'$ containing the vertex set $D$,
    \item[(2)] $G$ is bridgeless connected,
    \item[(3)] for every vertex $v\in V(G)\backslash D$ we have $\{v,f(v)\}\in E(G)$, where 
               $f:V(G')\backslash D\rightarrow D$ is the function from the first standard form of $(G',D)$, and
    \item[(4)] $G$ is vertex and edge-minimal with respect to properties (1), (2), and (3).
  \end{itemize}
\end{Definition}

\begin{Corollary}
  \label{cor_minimal_subgraph}
  If $G$ is a minimal subgraph of $(G',D)$ in first standard form, we have
  \begin{itemize}
    \item[(1)] $|V(G)|\le 5\cdot|D|-4$ and
    \item[(2)] there exists no chord $\{u,v\}\in E(G)$, where $\{u,v\}\cap D=\emptyset$.
  \end{itemize}
\end{Corollary}

\begin{Definition}
  Let $G$ be a minimal subgraph of $(G',D)$ in first standard form. We construct a graph $\tilde{G}$ from $G$
  by adding isolated triangles at vertices of $D$ such that $(\tilde{G},D)$ is in first standard form. We call
  $\tilde{G}$ a minimal completion and we say that $H$ is a minimal or an optimal orientation of $G$, if $H$
  is strongly connected and we have
  $$
    \overset{\longrightarrow}{diam}_{min}(\tilde{G})\ge \max\Big\{diam_0(H,D)+4,diam_1(H,D)+2,diam_2(H,D)\Big\}.
  $$
\end{Definition}

By considering the isolated triangles being associated to the vertices of the dominating set $D$ we can easily check, that every minimal subgraph $G$ of a pair $(G',D)$ in first standard form admits a minimal orientation $H$ and that we have $\overset{\longrightarrow}{diam}_{min}(G')\le \overset{\longrightarrow}{diam}_{min}\left(\tilde{G}\right)$. If $G$ does only fulfill
conditions (1)-(3) of Definition \ref{def_minimal_subgraph} then we may consider a minimal subgraph $G''$ of $(G',D)$, which
contains $G$ as a subgraph. With this we can call on orientation $H$ of $G$ minimal or optimal if it is minimal or optimal
for $G''$.

\begin{Definition}
  \label{def_critical}
  We call a pair $(G',D)$ in first standard form critical, if $\Xi(\gamma(G'))=\overset{\longrightarrow}{diam}_{min}\left(G'\right)$.
\end{Definition}

\begin{Definition}
  \label{def_critical_subgraph}
  We call a minimal subgraph $G$ of $(G',D)$ in first standard form critical if for a minimal orientation $H$ of $G$ we have
  $$
    \Xi(\gamma(G'))=\max\Big\{diam_0(H,D)+4,diam_1(H,D)+2,diam_2(H,D)\Big\}.
  $$
\end{Definition}

\noindent
Together with Lemma \ref{lemma_subgraph} we obtain:

\begin{Lemma}
  \begin{eqnarray*}
    \Xi(\gamma)&=& \max \Big\{ \min \Big\{ \,\max\left\{diam_0(H,D)+4,diam_1(H,D)+2,diam_2(H,D)\right\}\,:\,\\
               & & H\text{ is strongly connected orientation of }G\, \Big\} \,:\, G\text{ is critical minimal }\\
               & & \text{ subgraph of } (G',D) \text{ in first standard form with } |D|=\gamma\Big\}.
  \end{eqnarray*}
\end{Lemma}

Sometimes it is useful to know some basic facts about strongly connected orientations of graphs.

\begin{Lemma}
  \begin{itemize}
    \item[(1)] If $H$ is a strongly connected orientation of an undirected graph $G$ and $C$ is a directed cycle
               without repeated edges in $H$, then inverting of the edges of $C$ yields another strongly connected
               orientation of $G$.
    \item[(2)] If $H$ is a strongly connected orientation of an undirected graph $G$ and $P_1$ and $P_2$ are two
               edge-disjoint directed paths from $x$ to $y$, then inverting $P_2$ yields another strongly connected
               orientation of $G$.
    \item[(3)] If $H$ is a strongly connected orientation of an undirected graph $G$ then inverting all edges yields
               another strongly connected orientation with equal diameter.
  \end{itemize}
\end{Lemma}

\section{Reductions}
\label{sec_reductions}

In this section we will propose some reductions for critical minimal subgraphs $G$ of pairs $(G',D)$ in first standard form, in order to provide some tools for an inductive proof of a better upper bound on $\Xi(\gamma)$.

\begin{Lemma}
  \label{lemma_reduction_1}
  Let $G$ be a critical minimal subgraph of $(G',D)$ in first standard form with $\gamma=\gamma(G')=|D|\ge 3$. If $G$ contains
  vertices $x,y\in D$, $l_1,l_2,r_1,r_2\in V(G)\backslash D$, two edge disjoint paths $P_1=[x,l_1,r_1,y]$, $P_2=[x,l_2,r_2,y]$,
  all neighbors of $l_1,r_1$ are in $\{x,l_1,r_1,y\}$, and all neighbors of $l_2,r_2$ are in $\{x,l_2,r_2,y\}$, then we have
  $\Xi(\gamma)\le \Xi(\gamma-1)+3$.
\end{Lemma}
\begin{Proof}
  Let $\tilde{G}$ be the graph which arises from $G$ by deleting $l_1,l_2,r_1,r_2$ and identifying $x$ with $y$. Now let
  $\tilde{D}:=D\backslash\{y\}$ and $\tilde{H}$ be an arbitrary minimal orientation of $\tilde{G}$. Thus we have
  $diam_0\left(\tilde{H},\tilde{D}\right)\le\Xi(\gamma-1)+4$, $diam_1\left(\tilde{H},\tilde{D}\right)\le\Xi(\gamma-1)+2$, and
  $diam_2\left(\tilde{H},\tilde{D}\right)\le\Xi(\gamma-1)$. We construct an orientation $H$ of $G$ by directing the two paths
  $P_1$ and $P_2$
  in opposing  directions, and by taking the directions from $\tilde{H}$. Now we analyze the distance $d_H(u,v)$ in $H$ for all pairs
  $u,v\in V(G)$. If both $u$ and $v$ are in $\{l_1,l_2,r_1,r_2\}$, then we have $d_H(u,v)\le 5\le \Xi(\gamma-1)+3$. If none of $u$ and
  $v$ is in $\{l_1,l_2,r_1,r_2\}$, then we have $d_H(u,v)\le d_{\tilde{H}}(u,v)+3$. In the remaining case we have $d_H(u,v)\le
  d_{\tilde{H}}(u,v)+5$. Thus we have
  \begin{eqnarray*}
    diam_2(H,D) &\le&
    \max\Big\{diam_2\left(\tilde{H},\tilde{D}\right)+3,diam_1\left(\tilde{H},\tilde{D}\right)+5,5\Big\}\le\Xi(\gamma-1)+3,\\
    diam_1(H,D) &\le& \max\Big\{diam_1\left(\tilde{H},\tilde{D}\right)
    +3,diam_0\left(\tilde{H},\tilde{D}\right)+5,5\Big\}\le\Xi(\gamma-1)+1,\text{ and}\\
    diam_0(H,D) &\le& diam_0\left(\tilde{H},\tilde{D}\right)+3\le\Xi(\gamma-1)-1,
  \end{eqnarray*}
  which yields $\Xi(\gamma)\le\Xi(\gamma-1)+3$.
\end{Proof}

We remark that Lemma \ref{lemma_reduction_1} corresponds to a graph containing the left graph of Figure \ref{fig_gamma_2} as an induced subgraph, where the vertices corresponding to the empty circles have no further neighbors in the whole graph.

\begin{Lemma}
  \label{lemma_reduction_2}
  Let $G$ be a critical minimal subgraph of $(G',D)$ in first standard form with $\gamma=\gamma(G')=|D|\ge 3$. If $G$ 
  contains vertices $x,y,z\in D$, four edge disjoint paths $P_1=[x,v_1,v_2,v_3,y]$, $P_2=[y,v_4,v_5,v_6,z]$,
  $P_3=[x,u_1,u_2,y]$, $P_4=[y,u_3,u_4,z]$, and all edges being adjacent to vertices in
  $I:=\{v_1,v_2,v_3,v_4,v_5,v_6,u_1,u_2,u_3,u_4\}$ are contained in $P:=P_1\cup P_2\cup P_3\cup P_4$, then we have
  $\Xi(\gamma)\le \Xi(\gamma-2)+7$.
\end{Lemma}
\begin{Proof}
  At first we want to determine some structure information on the vertices $v_i$, $u_j$ and the adjacent edges. We have
  $f(v_1)=f(u_1)=x$, $f(v_3)=f(v_4)=f(u_2)=f(u_3)=y$, and $f(v_6)=f(u_4)=z$. Since all edges being adjacent to vertices in
  $I$ are contained in $P$ we have $f(v_2),f(v_5)\in\{x,y,z\}$. Some vertices may have several labels. By $v_i\sim$ we denote
  the set of labels which correspond to the same vertex as $v_i$. Similarly we define $u_i\sim$.

  Let us at first assume $|I|=10$, meaning, that each vertex has a unique label. In this case we may consider the edge
  $\{v_2,f(v_2)\}$ which is not contained in $P$ to see that $G$ would not be a minimal subgraph of $(G',D)$ in first standard
  form.

  Due to the $14$ pairwise different edges of $P$ and the information on the values of $f$ we have
  \begin{itemize}
    \item[(a)] $v_1\sim\subseteq\{v_1,v_5\}$, $v_3\sim\subseteq\{v_3,v_5\}$,
               $v_4\sim\subseteq\{v_2,v_4\}$, $v_6\sim\subseteq\{v_2,v_6\}$,
    \item[(b)] $u_1\sim\subseteq\{u_1,v_2,v_5\}$, $u_2\sim\subseteq\{u_2,v_2,v_5\}$,
               $u_3\sim\subseteq\{u_3,v_2,v_5\}$, $u_4\sim\subseteq\{u_4,v_2,v_5\}$,
    \item[(c)] $v_2\sim\subseteq\{v_2,v_4,v_5,v_6,u_1,u_2,u_3,u_4\}$, $v_5\sim\subseteq\{v_1,v_2,v_3,v_5,u_1,u_2,u_3,u_4\}$.
  \end{itemize}

  Next we assume $|I|=9$ which means that exactly one vertex in $I$ has two different labels and all other vertices have unique
  labels.
  \begin{itemize}
    \item[(1)] If $v_1=v_5$ then $v_2$, $v_3$, and $v_4$ could be deleted.
    \item[(2)] If $v_3=v_5$ then $v_4$ could be deleted.
    \item[(3)] If $u_1=v_2$ then by considering the edge $\{v_5,f(v_5)\}\notin P$ we could conclude that either $v_4$ or $v_6$
               could be deleted.
    \item[(4)] If $u_1=v_5$ then $u_2$ could be deleted.
    \item[(5)] If $u_2=v_2$ then by considering the edge $\{v_5,f(v_5)\}\notin P$ we could conclude that either $v_4$ or $v_6$
               could be deleted.
    \item[(6)] If $u_2=v_5$ then $v_4$ could be deleted.
    \item[(7)] If $v_2=v_5$ then $v_3$ and $v_4$ could be deleted.
  \end{itemize}
  Thus the vertices $v_1$, $v_3$, $u_1$, $u_2$ are unique. Using symmetry we conclude that also the vertices $v_4$, $v_6$, $u_3$, and
  $u_4$ are unique. Since we have also dealt with the only left possibility $v_2=v_5$ we can conclude $|I|\le 8$.

  We proceed similar as in the proof of Lemma \ref{lemma_reduction_1} and let $\tilde{G}$ be the graph arising from $G$ by deleting
  the vertices $u_i$, $v_i$, $y$ and by identifying $x$ and $z$. Obviously $\tilde{G}$ is bridgeless connected.
  Now let $\tilde{D}:=D\backslash\{y,z\}$ and $\tilde{H}$ be an
  arbitrary minimal orientation of $\tilde{G}$. Thus we have $diam_0\left(\tilde{H},\tilde{D}\right)\le\Xi(\gamma-2)-4$,
  $diam_1\left(\tilde{H},\tilde{D}\right)\le\Xi(\gamma-2)-2$, and $diam_2\left(\tilde{H},\tilde{D}\right)\le\Xi(\gamma-2)$.

  We construct an orientation $H$ of $G$ by directing the two pairs of paths $(P_1,P_3)$, $(P_2,P_4)$ in opposing  directions such
  that the arcs $[v_3,y]$, $[y,v_4]$ are directed different, by taking the directions from $\tilde{H}$ and by directing remaining
  edges arbitrarily.

  Now we analyze the distance $d_H(u,v)$ in $H$ for all pairs $u,v\in V(G)$. Due to $d_H(x,z),d_H(z,x)\le 7$, 
  $d_H(y,x),d_H(y,z),d_H(x,y),d_H(z,y)\le 4$ we have
  $d_H(u,v)\le d_{\tilde{H}}(u,v)+7$ for $u,v\notin I$. Now we consider $d_H(u,v)$ for $u,v\in I\cup\{x,y,z\}$. Due to
  $L:=|I\cup\{x,y,z\}|\le 11$ we clearly have $d_H(u,v)\le 10$. We assume $L=11$ since otherwise we would have $d_H(u,v)\le 9$.
  Now we have a closer look at the directed cycle $C:=P_1\circ P_4\circ P_2\circ P_3$ of length $14$ consisting of $11$ vertices.
  It is not possible to visit all $11$ vertices going along edges of the cycle $C$ without visiting a vertex twice. Thus we have
  $d_H(u,v)\le 9$ for $u,v\in I\cup\{x,y,z\}$. Summarizing our results gives
  \begin{eqnarray*}
    diam_2(H,D) &\le& \max\Big\{diam_2\left(\tilde{H},\tilde{D}\right)+7,diam_1\left(\tilde{H},\tilde{D}\right)+9,
    9\Big\}\le\Xi(\gamma-2)+7,\\
    diam_1(H,D) &\le& \max\Big\{diam_1\left(\tilde{H},\tilde{D}\right)+7,
    diam_0\left(\tilde{H},\tilde{D}\right)+9,9\Big\}\le\Xi(\gamma-2)+5,\text{ and}\\
    diam_0(H,D) &\le& diam_0\left(\tilde{H},\tilde{D}\right)+7\le\Xi(\gamma-2)+3,
  \end{eqnarray*}
  which yields $\Xi(\gamma)\le\Xi(\gamma-2)+7$.
\end{Proof}

We remark that Lemma \ref{lemma_reduction_2} corresponds to a graph containing the right graph of Figure \ref{fig_gamma_2} two times as an induced subgraph for $x,y,z\in D$ corresponding to the black circle, where the vertices corresponding to the empty circles have no further neighbors in the whole graph.

\begin{Lemma}
  \label{lemma_reduction_3}
  Let $G$ be a critical minimal subgraph of $(G',D)$ in first standard form with $\gamma=\gamma(G')=|D|\ge 3$ and $x$ a vertex
  contained in the dominating set $D$. If removing $x$ produces at least three connectivity components $C_1$, $C_2$, $C_3$,
  $\dots$, then we have 
  $$
    \Xi(\gamma)\le\max\Big\{\Xi(\gamma-i)+\Xi(i)-4\,:\, 1\le i\le \gamma-1\Big\}.
  $$
\end{Lemma}
\begin{Proof}
  Let $\tilde{C}_i$ be the induced subgraphs of $V(C_i)\cup\{x\}$ in $G$. We set $D_i=\{x\}\cup \left(V(C_i)\cap D\right)$ and
  $\gamma_i:=|D_i|-1$ so that we have $1+\sum_i \gamma_i=\gamma$. Since $G$ is a minimal subgraph we have $\gamma_i\ge 1$
  for all $i$. Now we choose arbitrary
  minimal orientations $\tilde{H}_i$ of the $\tilde{C}_i$. Thus we have $diam_0\left(\tilde{H}_i,D_i\right)\le\Xi(\gamma_i+1)-4$,
  $diam_1\left(\tilde{H}_i,D_i\right)\le\Xi(\gamma_i+1)-2$, and $diam_2\left(\tilde{H}_i,D_i\right)\le\Xi(\gamma_i+1)$ for all
  $i$. Since
  $\tilde{C}_i$ and $\tilde{C}_j$ are edge-disjoint for $i\neq j$ we can construct an orientation $H$ of $G$ by taking the
  directions of the $\tilde{H}_i$. Now we analyze the distance $d_H(u,v)$ in $H$ for all pairs $u,v\in V(G)$. If $u$ and $v$ are
  contained in the same component $\tilde{C}_i$ we have $d_H(u,v)=d_{\tilde{H}_i}(u,v)$. If $u$ is contained in $\tilde{C}_i$ and
  $v$ is contained in $\tilde{C}_j$, then we have $d_H(u,v)\le d_{\tilde{H}_i}(u,x)+d_{\tilde{H}_j}(x,v)$. Thus we have
  \begin{eqnarray*}
    diam_2(H,D) &\le& \max\Big\{diam_2\left(\tilde{H}_i,D_i\right),diam_1\left(\tilde{H}_i,D_i\right)
                +diam_1\left(\tilde{H}_j,D_j\right)\,:\,i\neq j\Big\}\\
                &\le& \max\Big\{\Xi(\gamma_i+1),\Xi(\gamma_i+1)+\Xi(\gamma_j+1)-4\,:\,i\neq j\Big\}\\
    diam_1(H,D) &\le& \max\Big\{diam_1\left(\tilde{H}_i,D_i\right),diam_1\left(\tilde{H}_i,D_i\right)
                +diam_0\left(\tilde{H}_j,D_j\right)\,:\, i\neq j\Big\}\\
                &\le& \max\Big\{\Xi(\gamma_i+1)-2,\Xi(\gamma_i+1)+\Xi(\gamma_j+1)-6\,:\,i\neq j\Big\},\text{ and}\\
    diam_0(H,D) &\le& \max\Big\{diam_0\left(\tilde{H}_i,D_i\right)+diam_0\left(\tilde{H}_j,D_j\right)\,:\, i\neq j\Big\}\\
                &\le& \max\Big\{\Xi(\gamma_i+1)+\Xi(\gamma_j+1)-8\,:\,i\neq j\Big\}.
  \end{eqnarray*}
  Since we have at least three connectivity components it holds $\gamma_i+\gamma_j\le \gamma -2$ for
  all $i\neq j$. Using this and $\Xi(n-1)\le\Xi(n)$ we conclude $\Xi(\gamma)\le\max\Big\{\Xi(\gamma-i)+\Xi(i)-4\,:\, 1\le i
  \le \gamma-1\Big\}$.
\end{Proof}

\begin{Lemma}
  \label{lemma_reduction_5}
  Let $G$ be a critical minimal subgraph of $(G',D)$ in first standard form with $\gamma=\gamma(G')=|D|\ge 3$ and $x$ a vertex
  not contained in the dominating set $D$. If removing $x$ produces at least three connectivity components $C_1$, $C_2$, $C_3$,
  $\dots$, then we have
  $$
    \Xi(\gamma)\le\max\Big\{\Xi(i)+\Xi(\gamma+1-i)-7,\Xi(i-1)+\Xi(\gamma+1-i)-4\,:\,2\le i\le\gamma-1\Big\}.
  $$
\end{Lemma}
\begin{Proof}
  W.l.o.g. let $f(x)$ be contained in $C_1$. Let $\tilde{C}_1$ be the induced subgraph of $V(C_1)\cup\{x\}$ in $G$ and
  $D_1=D\cap V(C_1)$.
  For $i\ge 2$ let $\tilde{C}_i$ be the induced subgraph of
  $V(C_i)\cup\{x\}$ in $G$ with additional vertices $y_i$, $z_i$, additional edges $\{x,y_i\}$, $\{x,z_i\}$,
  $\{y_i,z_i\}$, and $D_i=(V(C_i)\cap D)\cup \{z_i\}$. We set $\gamma_1=|D_1|\ge 1$ and $\gamma_i=|D_i|-1\ge 1$ for $i\ge 2$ so
  that we have $\sum_i \gamma_i=\gamma$. By $\tilde{H}_i$ we denote an optimal orientation of $C_i$. W.l.o.g. we assume that
  in $\tilde{H}_1$ the edge $\{f(x),x\}$ is directed from $f(x)$
  to $x$ and that for $i\ge 2$ in $\tilde{H}_i$ the edges $\{x,y_i\}$, $\{x,z_i\}$, $\{y_i,z_i\}$ are directed from $x$ to $y_i$,
  from $y_i$ to $z_i$ and from $z_i$ to $x$. Due to the minimality of the orientations $\tilde{H}_i$
  we have $diam_0\left(\tilde{H}_1,D_1\right)\le\Xi(\gamma_1)-4$, $diam_1\left(\tilde{H}_1,D_1\right)\le\Xi(\gamma_1)-2$,
  $diam_2\left(\tilde{H}_1,D_1\right)\le\Xi(\gamma_1)$, and for $i\ge 2$ we have $diam_0\left(\tilde{H}_i,D_i\right)
  \le\Xi(\gamma_i+1)-4$, $diam_1\left(\tilde{H}_i,D_i\right)\le\Xi(\gamma_i+1)-2$, $diam_2\left(\tilde{H}_i,D_i\right)
  \le\Xi(\gamma_i+1)$.

  We construct an orientation $H$ of $G$ by taking the directions of the common edges with the $\tilde{H}_i$. Now we analyze
  the distance $d_H(u,v)$ in $H$ for all pairs $u,v\in V(G)$. We only have to consider the cases where $u$ and $v$ are in
  different connectivity components. Let us first assume $u\in\tilde{C}_i$, $v\in\tilde{C}_j$ with $i,j\ge 2$. We have
  $$
    d_H(u,v)\le d_{\tilde{H}_i}(u,x)+d_{\tilde{H}_j}(x,v)\le d_{\tilde{H}_i}(u,z_i)-2+d_{\tilde{H}_j}(z_j,v)-1,
  $$
  since every directed path from a vertex $u\in V(G)$ to $z_i$ in $\tilde{H}_i$ uses the arcs $[x,y_i]$, $[y_i,z_i]$,
  and every directed path from $z_j$ to a vertex $v\in V(G)$ in $\tilde{H}_j$ uses the arc $[z_j,x]$. Now let $u$ be in
  $\tilde{C}_1$ and $v$ be in $\tilde{C}_i$ with $i\ge 2$. Since the edge $\{f(x),x\}$ is directed from $f(x)$ to $x$, both
  in $H$ and in $\tilde{H}_1$, we can conclude
  $$
    d_H(u,v)\le d_{\tilde{H}_1}(u,x)+d_{\tilde{H}_i}(x,v)\le d_{\tilde{H}_1}(u,f(x))+1+d_{\tilde{H}_i}(z_i,v)-1.
  $$
  If $u\in\tilde{C}_i$ with $i\ge 2$ and $v\in\tilde{C}_1$, then we similarly conclude
  $$
    d_H(u,v)\le d_{\tilde{H}_i}(u,x)+d_{\tilde{H}_1}(x,v)\le d_{\tilde{H}_i}(u,z_i)-2+d_{\tilde{H}_1}(x,v).
  $$

  \noindent
  Thus using $\Xi(i-1)\le \Xi(i)$ for $i\in\mathbb{N}$ and $\gamma_i+\gamma_j\le \gamma -1$ for
  all $i\neq j$ in total we have
  \begin{eqnarray*}
    diam_2(H,D) &\!\!\!\!\le\!\!\!\!& \max\Big\{diam_2\left(\tilde{H}_1,D_1\right),diam_2\left(\tilde{H}_i,D_i\right),
    diam_1\left(\tilde{H}_i,D_i\right)+diam_1\left(\tilde{H}_j,D_j\right)-3,\\
    && diam_1\left(\tilde{H}_1,D_1\right)+diam_1\left(\tilde{H}_i,D_i\right),diam_2\left(\tilde{H}_1,D_1\right)
       +diam_1\left(\tilde{H}_i,D_i\right)-2\Big\}\\
    &\!\!\!\!\le\!\!\!\!&
    \max\Big\{\Xi(\gamma-1),\Xi(\gamma_i+1)+\Xi(\gamma_j+1)-7,\Xi(\gamma_1)+\Xi(\gamma_i+1)-4\,:\,2\le i<j\Big\}\\
    &\!\!\!\!\le\!\!\!\!& \max\Big\{\Xi(i)+\Xi(\gamma+1-i)-7,\Xi(i-1)+\Xi(\gamma+1-i)-4\,:\,2\le i\le\gamma-1\Big\}\\
    diam_1(H,D) &\!\!\!\!\le\!\!\!\!& \max\Big\{diam_1\left(\tilde{H}_1,D_1\right),diam_1\left(\tilde{H}_i,D_i\right),
    diam_0\left(\tilde{H}_i,D_i\right)+diam_1\left(\tilde{H}_j,D_j\right)-3,\\
    && diam_0\left(\tilde{H}_1,D_1\right)+diam_1\left(\tilde{H}_i,D_i\right),diam_1\left(\tilde{H}_1,D_1\right)
       +diam_0\left(\tilde{H}_i,D_i\right),\\
    && diam_2\left(\tilde{H}_1,D_1\right)+diam_0\left(\tilde{H}_i,D_i\right)-2,
       diam_1\left(\tilde{H}_1,D_1\right)+diam_1\left(\tilde{H}_i,D_i\right)-2\Big\}\\
    &\!\!\!\!\le\!\!\!\!&
    \max\Big\{\Xi(\gamma-1)-2,\Xi(\gamma_i+1)+\Xi(\gamma_j+1)-9,\Xi(\gamma_1)+\Xi(\gamma_i+1)-6\,:\,2\le i<j\Big\}\\
    &\!\!\!\!\le\!\!\!\!& \max\Big\{\Xi(i)+\Xi(\gamma+1-i)-9,\Xi(i-1)+\Xi(\gamma+1-i)-6\,:\,2\le i\le\gamma-1\Big\}\\
    diam_0(H,D) &\!\!\!\!\le\!\!\!\!& \max\Big\{diam_0\left(\tilde{H}_1,D_1\right),diam_0\left(\tilde{H}_i,D_i\right),
    diam_0\left(\tilde{H}_i,D_i\right)+diam_0\left(\tilde{H}_j,D_j\right)-3,\\
    && diam_0\left(\tilde{H}_1,D_1\right)+diam_0\left(\tilde{H}_i,D_i\right),
       diam_1\left(\tilde{H}_1,D_1\right)+diam_0\left(\tilde{H}_i,D_i\right)-2\Big\}\\
    &\!\!\!\!\le\!\!\!\!&
    \max\Big\{\Xi(\gamma-1)-4,\Xi(\gamma_i+1)+\Xi(\gamma_j+1)-11,\Xi(\gamma_1)+\Xi(\gamma_i+1)-8\,:\,2\le i<j\Big\}\\
    &\!\!\!\!\le\!\!\!\!& \max\Big\{\Xi(i)+\Xi(\gamma+1-i)-11,\Xi(i-1)+\Xi(\gamma+1-i)-8\,:\,2\le i\le\gamma-1\Big\},
  \end{eqnarray*}
  which yields $\Xi(\gamma)\le\max\Big\{\Xi(i)+\Xi(\gamma+1-i)-7,\Xi(i-1)+\Xi(\gamma+1-i)-4\,:\,2\le i\le\gamma-1\Big\}$.
\end{Proof}

\begin{figure}[htp]
\begin{center}
  \setlength{\unitlength}{1cm}
  \begin{picture}(2.4,3.4)
    \put(0.2,0.2){\circle*{0.4}}
    \put(1.2,0.2){\circle{0.4}}
    \put(1.2,3.2){\circle{0.4}}
    \put(2.2,0.2){\circle{0.4}}
    \put(0.2,1.2){\circle{0.4}}
    \put(1.2,1.2){\circle{0.4}}
    \put(0.2,3.2){\circle{0.4}}
    \put(1.2,2.2){\circle*{0.4}}
    \put(2.2,3.2){\circle*{0.4}}
    \put(0.2,0.4){\vector(0,1){0.6}}
    \put(0.2,1.4){\vector(0,1){1.6}}
    \put(0.4,3.2){\vector(1,0){0.6}}
    \put(1.4,3.2){\vector(1,0){0.6}}
    \put(2.2,3.0){\vector(0,-1){2.6}}
    \put(2.0,0.2){\vector(-1,0){0.6}}
    \put(1.0,0.2){\vector(-1,0){0.6}}
    \put(1.2,0.4){\vector(0,1){0.6}}
    \put(1.2,1.4){\vector(0,1){0.6}}
    \put(1.065,2.335){\vector(-1,1){0.730}}
  \end{picture} 
  \quad
  \begin{picture}(2.4,3.4)
    \put(0.2,0.2){\circle*{0.4}}
    \put(0.2,1.2){\circle{0.4}}
    \put(0.2,2.2){\circle{0.4}}
    \put(0.2,3.2){\circle*{0.4}}
    \put(1.2,3.2){\circle{0.4}}
    \put(2.2,2.2){\circle{0.4}}
    \put(2.2,1.2){\circle*{0.4}}
    \put(2.2,0.2){\circle{0.4}}
    \put(1.2,0.2){\circle{0.4}}
    \put(0.2,0.4){\vector(0,1){0.6}}
    \put(0.2,1.4){\vector(0,1){0.6}}
    \put(0.2,2.4){\vector(0,1){0.6}}
    \put(0.4,3.2){\vector(1,0){0.6}}
    \put(1.335,3.065){\vector(1,-1){0.730}}
    \put(2.2,2.0){\vector(0,-1){0.6}}
    \put(2.2,1.0){\vector(0,-1){0.6}}
    \put(2.0,0.2){\vector(-1,0){0.6}}
    \put(1.0,0.2){\vector(-1,0){0.6}}
  \end{picture} 
  \quad
  \begin{picture}(3.4,2.4)
    \put(0.2,0.2){\circle{0.4}}
    \put(0.2,1.2){\circle{0.4}}
    \put(0.2,2.2){\circle*{0.4}}
    \put(1.2,2.2){\circle{0.4}}
    \put(1.2,1.2){\circle{0.4}}
    \put(1.2,0.2){\circle*{0.4}}
    \put(2.2,0.2){\circle{0.4}}
    \put(2.2,2.2){\circle{0.4}}
    \put(3.2,1.2){\circle*{0.4}}
    \put(0.2,0.4){\vector(0,1){0.6}}
    \put(0.2,1.4){\vector(0,1){0.6}}
    \put(0.4,2.2){\vector(1,0){0.6}}
    \put(1.2,2.0){\vector(0,-1){0.6}}
    \put(1.2,1.0){\vector(0,-1){0.6}}
    \put(1.0,0.2){\vector(-1,0){0.6}}
    \put(1.335,1.335){\vector(1,1){0.730}}
    \put(2.335,2.065){\vector(1,-1){0.730}}
    \put(3.065,1.065){\vector(-1,-1){0.730}}
    \put(2.065,0.335){\vector(-1,1){0.730}}
  \end{picture} 
  \quad
  \begin{picture}(2.4,2.4)
    \put(0.2,0.2){\circle*{0.4}}
    \put(0.2,1.2){\circle{0.4}}
    \put(0.2,2.2){\circle{0.4}}
    \put(1.2,0.2){\circle{0.4}}
    \put(1.2,1.2){\circle*{0.4}}
    \put(1.2,2.2){\circle{0.4}}
    \put(2.2,0.2){\circle{0.4}}
    \put(2.2,1.2){\circle{0.4}}
    \put(2.2,2.2){\circle*{0.4}}
    \put(0.2,0.4){\vector(0,1){0.6}}
    \put(0.2,1.4){\vector(0,1){0.6}}
    \put(0.4,2.2){\vector(1,0){0.6}}
    \put(1.4,2.2){\vector(1,0){0.6}}
    \put(2.2,2.0){\vector(0,-1){0.6}}
    \put(2.2,1.0){\vector(0,-1){0.6}}
    \put(2.0,0.2){\vector(-1,0){0.6}}
    \put(1.0,0.2){\vector(-1,0){0.6}}
    \put(0.335,2.065){\vector(1,-1){0.73}}
    \put(1.335,1.065){\vector(1,-1){0.73}}
  \end{picture}\\[4mm] 
  \begin{picture}(2.4,3.4)
    \put(0.2,0.2){\circle*{0.4}}
    \put(0.2,1.2){\circle{0.4}}
    \put(0.2,2.2){\circle{0.4}}
    \put(0.2,3.2){\circle{0.4}}
    \put(1.2,3.2){\circle*{0.4}}
    \put(1.2,2.2){\circle{0.4}}
    \put(2.2,2.2){\circle*{0.4}}
    \put(2.2,0.2){\circle{0.4}}
    \put(1.2,0.2){\circle{0.4}}
    \put(0.2,0.4){\vector(0,1){0.6}}
    \put(0.2,1.4){\vector(0,1){0.6}}
    \put(0.4,2.2){\vector(1,0){0.6}}
    \put(1.4,2.2){\vector(1,0){0.6}}
    \put(2.2,2.0){\vector(0,-1){1.6}}
    \put(2.0,0.2){\vector(-1,0){0.6}}
    \put(1.0,0.2){\vector(-1,0){0.6}}
    \put(0.335,2.335){\vector(1,1){0.73}}
    \put(1.0,3.2){\vector(-1,0){0.6}}
    \put(0.2,3.0){\vector(0,-1){0.6}}
  \end{picture} 
  \quad
  \begin{picture}(3.4,3.4)
    \put(0.2,0.2){\circle*{0.4}}
    \put(0.2,2.2){\circle{0.4}}
    \put(1.2,2.2){\circle{0.4}}
    \put(1.2,1.2){\circle*{0.4}}
    \put(2.2,1.2){\circle{0.4}}
    \put(3.2,0.2){\circle{0.4}}
    \put(3.2,2.2){\circle{0.4}}
    \put(3.2,3.2){\circle{0.4}}
    \put(2.2,3.2){\circle*{0.4}}
    \put(0.2,0.4){\vector(0,1){1.6}}
    \put(0.4,2.2){\vector(1,0){0.6}}
    \put(1.4,2.2){\vector(1,0){1.6}}
    \put(3.2,2.0){\vector(0,-1){1.6}}
    \put(3.0,0.2){\vector(-1,0){2.6}}
    \put(3.065,2.335){\vector(-1,1){0.73}}
    \put(2.4,3.2){\vector(1,0){0.6}}
    \put(3.2,3.0){\vector(0,-1){0.6}}
    \put(3.065,0.335){\vector(-1,1){0.73}}
    \put(2.0,1.2){\vector(-1,0){0.6}}
    \put(1.2,1.4){\vector(0,1){0.6}}
  \end{picture} 
  \quad
  \begin{picture}(3.4,2.4)
    \put(1.2,0.2){\circle{0.4}}
    \put(1.2,1.2){\circle{0.4}}
    \put(0.2,1.2){\circle{0.4}}
    \put(0.2,2.2){\circle*{0.4}}
    \put(2.2,0.2){\circle*{0.4}}
    \put(2.2,1.2){\circle{0.4}}
    \put(2.2,2.2){\circle{0.4}}
    \put(3.2,1.2){\circle{0.4}}
    \put(3.2,2.2){\circle*{0.4}}
    \put(1.2,0.4){\vector(0,1){0.6}}
    \put(0.4,1.2){\vector(1,0){0.6}}
    \put(1.4,1.2){\vector(1,0){0.6}}
    \put(2.2,1.0){\vector(0,-1){0.6}}
    \put(2.0,0.2){\vector(-1,0){0.6}}
    \put(2.2,1.4){\vector(0,1){0.6}}
    \put(2.4,2.2){\vector(1,0){0.6}}
    \put(3.2,2.0){\vector(0,-1){0.6}}
    \put(3.0,1.2){\vector(-1,0){0.6}}
    \put(1.065,1.335){\vector(-1,1){0.73}}
    \put(0.2,2.0){\vector(0,-1){0.6}}
  \end{picture} 
  \quad
  \begin{picture}(2.4,3.4)
    \put(0.2,0.2){\circle*{0.4}}
    \put(0.2,2.2){\circle{0.4}}
    \put(1.2,2.2){\circle{0.4}}
    \put(1.2,1.2){\circle*{0.4}}
    \put(2.2,1.2){\circle{0.4}}
    \put(2.2,0.2){\circle{0.4}}
    \put(2.2,2.2){\circle{0.4}}
    \put(2.2,3.2){\circle{0.4}}
    \put(1.2,3.2){\circle*{0.4}}
    \put(0.2,0.4){\vector(0,1){1.6}}
    \put(0.4,2.2){\vector(1,0){0.6}}
    \put(1.4,2.2){\vector(1,0){0.6}}
    \put(2.0,0.2){\vector(-1,0){1.6}}
    \put(2.065,2.335){\vector(-1,1){0.73}}
    \put(1.4,3.2){\vector(1,0){0.6}}
    \put(2.2,3.0){\vector(0,-1){0.6}}
    \put(2.0,1.2){\vector(-1,0){0.6}}
    \put(1.2,1.4){\vector(0,1){0.6}}
    \put(2.2,2.0){\vector(0,-1){0.6}}
    \put(2.2,1.0){\vector(0,-1){0.6}}
  \end{picture}\\[4mm] 
  \begin{picture}(4.4,1.4)
    \put(0.2,0.2){\circle{0.4}}
    \put(1.2,0.2){\circle{0.4}}
    \put(2.2,0.2){\circle*{0.4}}
    \put(3.2,0.2){\circle{0.4}}
    \put(4.2,0.2){\circle{0.4}}
    \put(0.2,1.2){\circle*{0.4}}
    \put(1.2,1.2){\circle{0.4}}
    \put(3.2,1.2){\circle{0.4}}
    \put(4.2,1.2){\circle*{0.4}}
    \put(1.2,0.4){\vector(0,1){0.6}}
    \put(1.4,1.2){\vector(1,0){1.6}}
    \put(3.2,1.0){\vector(0,-1){0.6}}
    \put(3.0,0.2){\vector(-1,0){0.6}}
    \put(2.0,0.2){\vector(-1,0){0.6}}
    \put(0.4,1.2){\vector(1,0){0.6}}
    \put(0.2,0.4){\vector(0,1){0.6}}
    \put(1.065,1.065){\vector(-1,-1){0.73}}
    \put(3.4,1.2){\vector(1,0){0.6}}
    \put(4.2,1.0){\vector(0,-1){0.6}}
    \put(4.065,0.335){\vector(-1,1){0.73}}
  \end{picture} 
  \quad
  \begin{picture}(4.4,1.4)
    \put(0.2,0.2){\circle*{0.4}}
    \put(1.2,0.2){\circle{0.4}}
    \put(2.2,0.2){\circle*{0.4}}
    \put(3.2,0.2){\circle{0.4}}
    \put(4.2,0.2){\circle*{0.4}}
    \put(0.2,1.2){\circle{0.4}}
    \put(1.2,1.2){\circle{0.4}}
    \put(3.2,1.2){\circle{0.4}}
    \put(4.2,1.2){\circle{0.4}}
    \put(0.2,0.4){\vector(0,1){0.6}}
    \put(0.4,1.2){\vector(1,0){0.6}}
    \put(1.4,1.2){\vector(1,0){1.6}}
    \put(3.4,1.2){\vector(1,0){0.6}}
    \put(4.2,1.0){\vector(0,-1){0.6}}
    \put(4.0,0.2){\vector(-1,0){0.6}}
    \put(3.2,0.4){\vector(0,1){0.6}}
    \put(1.2,1.0){\vector(0,-1){0.6}}
    \put(1.0,0.2){\vector(-1,0){0.6}}
    \put(3.065,1.065){\vector(-1,-1){0.73}}
    \put(2.065,0.335){\vector(-1,1){0.73}}
  \end{picture}\\[4mm] 
  \begin{picture}(4.4,2.4)
    \put(2.2,0.2){\circle{0.4}}
    \put(2.2,1.2){\circle*{0.4}}
    \put(1.2,1.2){\circle{0.4}}
    \put(3.2,1.2){\circle{0.4}}
    \put(2.2,2.2){\circle{0.4}}
    \put(1.2,2.2){\circle*{0.4}}
    \put(3.2,2.2){\circle*{0.4}}
    \put(0.2,2.2){\circle{0.4}}
    \put(4.2,2.2){\circle{0.4}}
    \put(2.2,0.4){\vector(0,1){0.6}}
    \put(2.2,1.4){\vector(0,1){0.6}}
    \put(2.065,0.335){\vector(-1,1){0.73}}
    \put(1.065,1.335){\vector(-1,1){0.73}}
    \put(4.065,2.065){\vector(-1,-1){0.73}}
    \put(3.065,1.065){\vector(-1,-1){0.73}}
    \put(0.4,2.2){\vector(1,0){0.6}}
    \put(3.4,2.2){\vector(1,0){0.6}}
    \put(1.335,1.335){\vector(1,1){0.73}}
    \put(2.335,2.065){\vector(1,-1){0.73}}
    \put(1.2,2.0){\vector(0,-1){0.6}}
    \put(3.2,1.4){\vector(0,1){0.6}}
  \end{picture} 
  \quad
  \begin{picture}(4.4,2.4)
    \put(1.2,0.2){\circle{0.4}}
    \put(2.2,0.2){\circle*{0.4}}
    \put(2.2,1.2){\circle{0.4}}
    \put(1.2,2.2){\circle{0.4}}
    \put(3.2,2.2){\circle{0.4}}
    \put(1.2,1.2){\circle*{0.4}}
    \put(0.2,1.2){\circle{0.4}}
    \put(3.2,1.2){\circle*{0.4}}
    \put(4.2,1.2){\circle{0.4}}
    \put(2.2,1.0){\vector(0,-1){0.6}}
    \put(2.0,0.2){\vector(-1,0){0.6}}
    \put(1.335,0.335){\vector(1,1){0.73}}
    \put(1.0,1.2){\vector(-1,0){0.6}}
    \put(1.2,2.0){\vector(0,-1){0.6}}
    \put(1.4,2.2){\vector(1,0){1.6}}
    \put(3.2,2.0){\vector(0,-1){0.6}}
    \put(3.4,1.2){\vector(1,0){0.6}}
    \put(0.335,1.335){\vector(1,1){0.73}}
    \put(2.065,1.335){\vector(-1,1){0.73}}
    \put(3.065,2.065){\vector(-1,-1){0.73}}
    \put(4.065,1.335){\vector(-1,1){0.73}}
  \end{picture}\\[4mm] 
  \begin{picture}(4.4,2.4)
    \put(1.2,0.2){\circle{0.4}}
    \put(2.2,0.2){\circle*{0.4}}
    \put(3.2,0.2){\circle{0.4}}
    \put(0.2,1.2){\circle{0.4}}
    \put(2.2,1.2){\circle{0.4}}
    \put(4.2,1.2){\circle{0.4}}
    \put(0.2,2.2){\circle*{0.4}}
    \put(1.2,2.2){\circle{0.4}}
    \put(3.2,2.2){\circle{0.4}}
    \put(4.2,2.2){\circle*{0.4}}
    \put(2.2,1.0){\vector(0,-1){0.6}}
    \put(2.0,0.2){\vector(-1,0){0.6}}
    \put(1.065,0.335){\vector(-1,1){0.73}}
    \put(0.2,1.4){\vector(0,1){0.6}}
    \put(0.4,2.2){\vector(1,0){0.6}}
    \put(1.335,2.065){\vector(1,-1){0.73}}
    \put(2.335,1.335){\vector(1,1){0.73}}
    \put(3.4,2.2){\vector(1,0){0.6}}
    \put(4.2,2.0){\vector(0,-1){0.6}}
    \put(4.065,1.065){\vector(-1,-1){0.73}}
    \put(3.0,0.2){\vector(-1,0){0.6}} 
  \end{picture}\quad 
  \begin{picture}(4.4,2.4)
    \put(1.2,0.2){\circle{0.4}}
    \put(2.2,0.2){\circle*{0.4}}
    \put(3.2,0.2){\circle{0.4}}
    \put(0.2,1.2){\circle{0.4}}
    \put(2.2,1.2){\circle{0.4}}
    \put(4.2,1.2){\circle{0.4}}
    \put(0.2,2.2){\circle*{0.4}}
    \put(1.2,2.2){\circle{0.4}}
    \put(3.2,1.2){\circle{0.4}}
    \put(4.2,2.2){\circle*{0.4}}
    \put(2.2,1.0){\vector(0,-1){0.6}}
    \put(2.0,0.2){\vector(-1,0){0.6}}
    \put(1.065,0.335){\vector(-1,1){0.73}}
    \put(0.2,1.4){\vector(0,1){0.6}}
    \put(0.4,2.2){\vector(1,0){0.6}}
    \put(1.335,2.065){\vector(1,-1){0.73}}
    \put(2.4,1.2){\vector(1,0){0.6}}
    \put(3.2,1.0){\vector(0,-1){0.6}}
    \put(3.335,1.335){\vector(1,1){0.73}}
    \put(4.2,2.0){\vector(0,-1){0.6}}
    \put(4.0,1.2){\vector(-1,0){0.6}}
    \put(3.0,0.2){\vector(-1,0){0.6}} 
  \end{picture}\quad 
  \begin{picture}(3.4,2.4)
    \put(0.2,0.2){\circle{0.4}}
    \put(0.2,1.2){\circle{0.4}}
    \put(0.2,2.2){\circle*{0.4}}
    \put(1.2,0.2){\circle*{0.4}}
    \put(1.2,1.2){\circle{0.4}}
    \put(1.2,2.2){\circle{0.4}}
    \put(2.2,0.2){\circle{0.4}}
    \put(2.2,1.2){\circle{0.4}}
    \put(3.2,0.2){\circle*{0.4}}
    \put(3.2,1.2){\circle{0.4}}
    \put(1.0,0.2){\vector(-1,0){0.6}}
    \put(0.2,0.4){\vector(0,1){0.6}}
    \put(0.2,1.4){\vector(0,1){0.6}}
    \put(0.4,2.2){\vector(1,0){0.6}}
    \put(1.2,2.0){\vector(0,-1){0.6}}
    \put(1.2,1.0){\vector(0,-1){0.6}}
    \put(1.4,1.2){\vector(1,0){0.6}}
    \put(2.4,1.2){\vector(1,0){0.6}}
    \put(3.2,1.0){\vector(0,-1){0.6}}
    \put(3.0,0.2){\vector(-1,0){0.6}}
    \put(2.2,0.4){\vector(0,1){0.6}}
    \put(2.065,1.065){\vector(-1,-1){0.73}}
  \end{picture}\\[4mm] 
  \begin{picture}(4.4,2.4)
    \put(1.2,0.2){\circle{0.4}}
    \put(1.2,1.2){\circle{0.4}}
    \put(1.2,2.2){\circle*{0.4}}
    \put(0.2,2.2){\circle{0.4}}
    \put(2.2,0.2){\circle*{0.4}}
    \put(2.2,1.2){\circle{0.4}}
    \put(3.2,0.2){\circle{0.4}}
    \put(3.2,1.2){\circle{0.4}}
    \put(3.2,2.2){\circle*{0.4}}
    \put(4.2,2.2){\circle{0.4}}
    \put(2.0,0.2){\vector(-1,0){0.6}}
    \put(1.2,0.4){\vector(0,1){0.6}}
    \put(1.2,1.4){\vector(0,1){0.6}}
    \put(1.0,2.2){\vector(-1,0){0.6}}
    \put(0.335,2.065){\vector(1,-1){0.73}}
    \put(1.4,1.2){\vector(1,0){0.6}}
    \put(2.2,1.0){\vector(0,-1){0.6}}
    \put(2.4,1.2){\vector(1,0){0.6}}
    \put(3.2,1.4){\vector(0,1){0.6}}
    \put(3.2,1.0){\vector(0,-1){0.6}}
    \put(3.0,0.2){\vector(-1,0){0.6}}
    \put(3.4,2.2){\vector(1,0){0.6}}
    \put(4.065,2.065){\vector(-1,-1){0.73}}
  \end{picture}\quad\quad 
  \begin{picture}(4.4,2.4)
    \put(1.2,0.2){\circle{0.4}}
    \put(0.2,0.2){\circle{0.4}}
    \put(0.2,1.2){\circle*{0.4}}
    \put(2.2,0.2){\circle*{0.4}}
    \put(2.2,1.2){\circle{0.4}}
    \put(1.2,1.2){\circle{0.4}}
    \put(3.2,0.2){\circle{0.4}}
    \put(3.2,1.2){\circle{0.4}}
    \put(3.2,2.2){\circle*{0.4}}
    \put(4.2,2.2){\circle{0.4}}
    \put(1.4,0.2){\vector(1,0){0.6}}
    \put(2.2,0.4){\vector(0,1){0.6}}
    \put(1.0,0.2){\vector(-1,0){0.6}}
    \put(0.2,0.4){\vector(0,1){0.6}}
    \put(0.4,1.2){\vector(1,0){0.6}}
    \put(1.2,1.0){\vector(0,-1){0.6}}
    \put(2.4,1.2){\vector(1,0){0.6}}
    \put(3.2,1.4){\vector(0,1){0.6}}
    \put(3.4,2.2){\vector(1,0){0.6}}
    \put(4.065,2.065){\vector(-1,-1){0.73}}
    \put(3.2,1.0){\vector(0,-1){0.6}}
    \put(3.0,0.2){\vector(-1,0){0.6}}
    \put(2.065,1.065){\vector(-1,-1){0.73}}
  \end{picture}
  \caption{The orientations for the proof of Lemma \ref{lemma_gamma_3}.}
  \label{fig_gamma_3}
\end{center}
\end{figure}
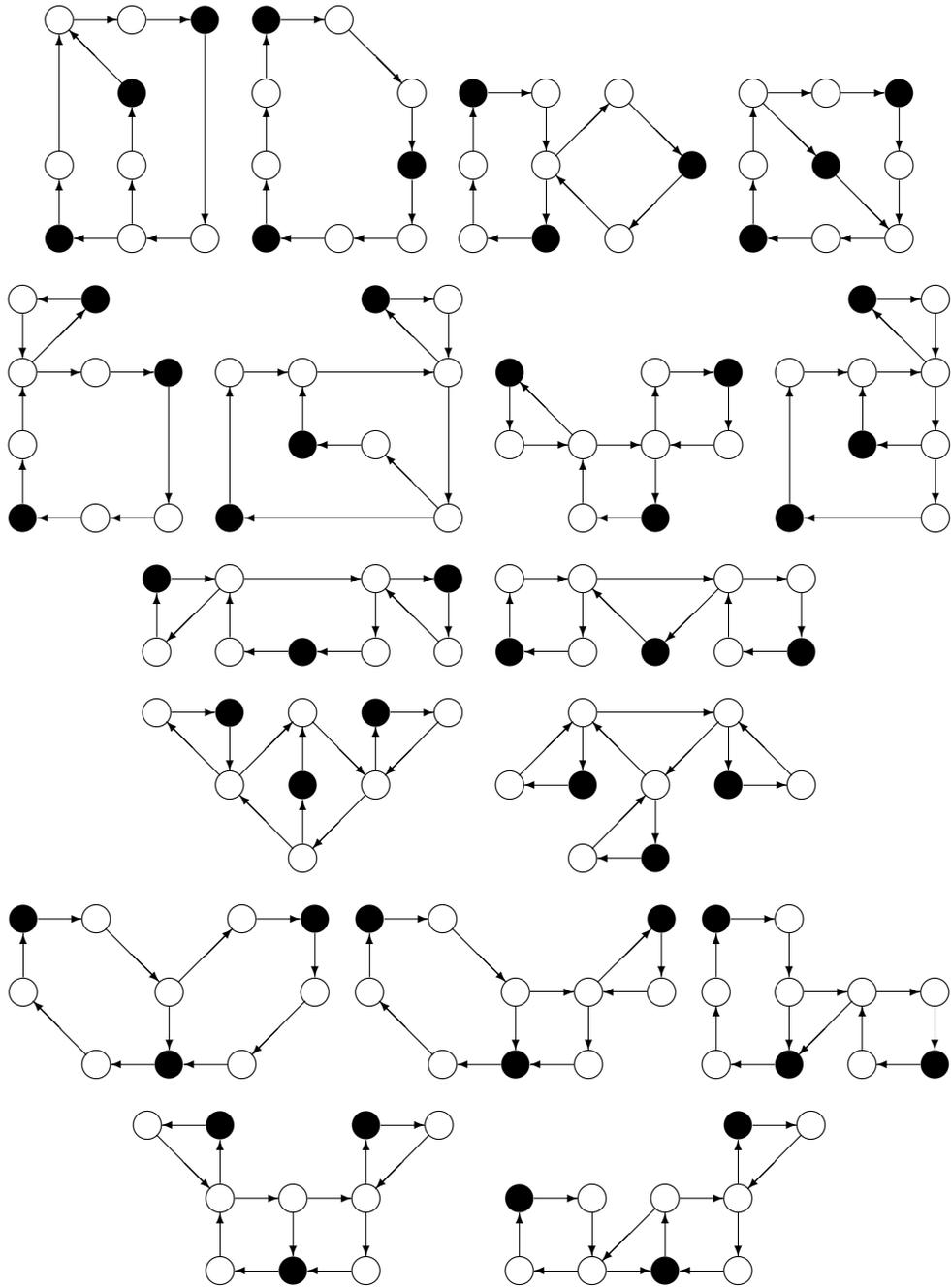

Now we are ready to determine the next exact value of $\Xi(\gamma)$:

\begin{Lemma}
  \label{lemma_gamma_3}
  $$
    \Xi(3)=11.
  $$
\end{Lemma}
\begin{Proof}
  The last example from Figure \ref{fig_the_bad_examples} gives $\Xi(3)\ge 11$. Going through the cases of the subgraph
  construction being described in front of Lemma \ref{lemma_small_exact_values} we are able to explicitly construct
  a finite list of possible subgraphs for $\gamma=3$. This fall differentiation is a bit laborious but not difficult. We can
  assume that these graphs $G$ are minimal subgraphs of a suitable pair $(G',D)$ in first standard form. During
  our construction we can drop all graphs which are not minimal, e.~g. graphs containing a chord where no end vertex lies
  in the dominating set $D$. Doing this we obtain a list of $24$ non-isomorphic minimal subgraphs. In Figure \ref{fig_gamma_3}
  we give suitable orientations for the cases, where we can not apply Lemma \ref{lemma_reduction_1}, Lemma \ref{lemma_reduction_2},
  or Lemma \ref{lemma_reduction_5}.
\end{Proof}
%

Going over the proofs of the previous lemmas again, we can conclude some further, in some sense weaker, reduction results.
Similarly as in Lemma \ref{lemma_reduction_2} we can prove:

\begin{Lemma}
  \label{lemma_weak_reduction_1}
  Let $G$ be a critical minimal subgraph of $(G',D)$ in first standard form with $\gamma=\gamma(G')=|D|\ge 3$. If $G$ 
  contains vertices $x,y,\in D$, two edge disjoint paths $P_1=[x,u_1,u_2,u_3,y]$,  $P_2=[x,v_1,v_2,y]$, and all edges 
  being adjacent to vertices in $I:=\{u_1,u_2,u_3,v_1,v_2\}$ are contained in $P_1\cup P_2$, then we have $\Xi(\gamma)
  \le \Xi(\gamma-1)+4$.
\end{Lemma}

\begin{Lemma}
  \label{lemma_weak_reduction_2}
  Let $G$ be a critical minimal subgraph of $(G',D)$ in first standard form with $\gamma=\gamma(G')=|D|\ge 3$ and $x$ a vertex
  contained in the dominating set $D$. If removing $x$ produces two connectivity components $C_1$ and $C_2$ then we have 
  $$
    \Xi(\gamma)\le\max\Big\{\Xi(\gamma+1-i)+\Xi(i)-4\,:\, 2\le i\le \gamma-1\Big\}.
  $$
\end{Lemma}
\begin{Proof}
  We can rephrase most of the proof of Lemma \ref{lemma_reduction_3}. Our estimations on $diam_i(H,D)$ remain valid. Since we
  only have two connectivity components we do not have $\gamma_i+\gamma_j\le\gamma-2$ for $i\neq j$. Instead we have
  $\gamma_1+\gamma_2=\gamma-1$ and $\gamma_1,\gamma_2\le\gamma-2$. Combining this with $\Xi(n-1)\le\Xi(n)$ we obtain the stated
  upper bound.
\end{Proof}

\begin{Lemma}
  \label{lemma_weak_reduction_3}
  Let $G$ be a critical minimal subgraph of $(G',D)$ in first standard form with $\gamma=\gamma(G')=|D|\ge 3$ and $x$ a vertex
  not contained in the dominating set $D$. If removing $x$ produces at least two connectivity components $C_1$, $C_2$ where
  $f(x)\in C_1$ and $|V(C_1)\cap D|\ge 2$ then we have
  $$
    \Xi(\gamma)\le\max\Big\{\Xi(i)+\Xi(\gamma+1-i)-4\,:\, 2\le i\le \gamma-1\Big\}.
  $$
\end{Lemma}
\begin{Proof}
  We can rephrase most of the proof of Lemma \ref{lemma_reduction_5}. Using $\Xi(i-1)\le\Xi(i)$ for all $i\in\mathbb{N}$ and
  the fact that we have exactly two connectivity components $C_1$ and $C_2$ yields
  \begin{eqnarray*}
    diam_2(H,D) &\!\!\!\!\le\!\!\!\!&
    \max\Big\{\Xi(\gamma-1),\Xi(\gamma_1)+\Xi(\gamma_2+1)-4\Big\}\\
    diam_1(H,D) &\!\!\!\!\le\!\!\!\!&
    \max\Big\{\Xi(\gamma-1)-2,\Xi(\gamma_1)+\Xi(\gamma_2+1)-6\Big\}\\
    diam_0(H,D) &\!\!\!\!\le\!\!\!\!&
    \max\Big\{\Xi(\gamma-1)-4,\Xi(\gamma_1)+\Xi(\gamma_2+1)-8\Big\}.
  \end{eqnarray*}
  Due to $\Xi(i-1)\le\Xi(i)$, $2\le y_1\le\gamma-1$, and $1\le\gamma_2\le\gamma-1$ we have
  $$
    \Xi(\gamma)\le\max\Big\{\Xi(i)+\Xi(\gamma+1-i)-4\,:\, 2\le i\le \gamma-1\Big\}.
  $$
\end{Proof}

We would like to remark that Lemmas \ref{lemma_reduction_1}, \ref{lemma_reduction_2}, \ref{lemma_reduction_3},
\ref{lemma_reduction_5} can be used in an induction proof of Conjecture \ref{main_conj}, whereas Lemmas \ref{lemma_weak_reduction_1}, \ref{lemma_weak_reduction_2}, \ref{lemma_weak_reduction_3} can only be used in an induction proof of Theorem \ref{thm_four_gamma}.

\begin{figure}[htp]
\begin{center}
  \setlength{\unitlength}{1cm}
  \begin{picture}(4.4,2.4)
    \put(0.2,1.2){\circle*{0.4}}
    \put(1.2,2.2){\circle{0.4}}
    \put(2.2,1.2){\circle{0.4}}
    \put(3.2,0.2){\circle{0.4}}
    \put(3.2,2.2){\circle{0.4}}
    \put(4.2,1.2){\circle*{0.4}}
    \put(-0.1,0.7){$f(x)$}
    \put(2.115,0.7){$x$}
    \put(4.115,0.7){$z$}
    \put(1.45,2.15){$w$}
    \put(3.45,2.15){$y_1$}
    \put(3.45,0.15){$y_2$}
    \put(0.2,1.2){\line(1,0){1.8}}
    \put(0.2,1.2){\line(1,1){0.865}}
    \put(2.065,1.335){\line(-1,1){0.730}}
    \put(3.335,2.065){\line(1,-1){0.750}}
    \put(3.335,0.335){\line(1,1){0.750}}
    \put(3.065,2.065){\line(-1,-1){0.730}}
    \put(3.065,0.335){\line(-1,1){0.730}}
  \end{picture}
  \quad\!
  \begin{picture}(2.4,2.4)
    \put(1.2,1.2){\circle{0.4}}
    \put(1.2,0.2){\circle{0.4}}
    \put(1.2,2.2){\circle{0.4}}
    \put(2.2,1.2){\circle*{0.4}}
    \put(0.75,1.15){$x$}
    \put(2.115,0.7){$z$}
    \put(1.45,2.15){$y_1$}
    \put(1.45,0.15){$y_2$}
    \put(2.2,1.2){\vector(-1,1){0.87}}
    \put(1.335,0.335){\line(1,1){0.750}}
    \put(1.2,1.4){\line(0,1){0.6}}
    \put(1.2,1.0){\line(0,-1){0.6}}
    \put(1.4,1.2){\line(1,0){0.8}}
  \end{picture}
  \quad\quad\,\,
  \begin{picture}(4.4,2.4)
    \put(0.2,1.2){\circle*{0.4}}
    \put(1.2,2.2){\circle{0.4}}
    \put(2.2,1.2){\circle{0.4}}
    \put(3.2,0.2){\circle{0.4}}
    \put(3.2,2.2){\circle{0.4}}
    \put(4.2,1.2){\circle*{0.4}}
    \put(-0.1,0.7){$f(x)$}
    \put(2.115,0.7){$x$}
    \put(4.115,0.7){$z$}
    \put(1.45,2.15){$w$}
    \put(3.45,2.15){$y_1$}
    \put(3.45,0.15){$y_2$}
    \put(0.2,1.2){\vector(1,0){1.8}}
    \put(1.065,2.065){\vector(-1,-1){0.730}}
    \put(2.065,1.335){\vector(-1,1){0.730}}
    \put(3.335,2.065){\vector(1,-1){0.730}}
    \put(4.065,1.065){\vector(-1,-1){0.730}}
    \put(2.335,1.335){\vector(1,1){0.730}}
    \put(3.065,0.335){\vector(-1,1){0.730}}
  \end{picture}
  \caption{The situation of Lemma \ref{lemma_weak_reduction_4} if we can not apply Lemma \ref{lemma_weak_reduction_3}.}
  \label{fig_situation_weak_reduction_4}
\end{center}
\end{figure}
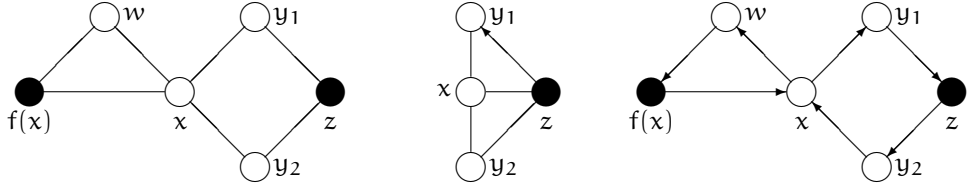

In order to prove Theorem \ref{thm_four_gamma} we need some further reduction Lemmas.

\begin{Lemma}
  \label{lemma_weak_reduction_4}
  Let $G$ be a critical minimal subgraph of $(G',D)$ in first standard form with $\gamma=\gamma(G')=|D|\ge 3$ and $x$ a vertex
  not contained in the dominating set $D$. If removing $x$ produces at least two connectivity components $C_1$, $C_2$, where
  $f(x)\in C_1$ and their exist $y_1\neq y_2\in V(G)\backslash D$ fulfilling $f(y_1)=f(y_2)$ and $\{x,y_1\},\{x,y_2\}\in E(G)$
  then we either can apply Lemma \ref{lemma_weak_reduction_3} or we have $\Xi(\gamma)\le\Xi(\gamma-1)+4$.
\end{Lemma}
\begin{Proof}
  If $|V(C_1)\cap D|\ge 2$ we can apply Lemma \ref{lemma_weak_reduction_3} thus we may assume $|V(C_1)\cap D|=1$. Since $G$ is
  a minimal subgraph, we have $V(C_1)=\{f(x),w\}$ and the neighbors of $f(x)$ and $w$ in $G$ are contained in $\{f(x),w,x\}$.
  As an abbreviation we set $f(y_1)=f(y_2)=z\in D$. See the left drawing in Figure \ref{fig_situation_weak_reduction_4} for a
  graphical representation of the situation. Now we consider the subgraph $\tilde{C}_2$ consisting of the induced subgraph of
  $V(C_2)\cup\{x\}$ with the additional edge $\{x,f(y_1)\}$. Let $H_2$ be an optimal orientation of $\tilde{C}_2$, where we
  assume that the arc $[z,y_1]$ is directed from $z$ to $y_1$, see the middle graph of Figure \ref{fig_situation_weak_reduction_4}.
  Now we construct an orientation $H$ of $G$ by taking the directions from $H_2$ and redirecting some edges. We direct $x$ to
  $w$, $w$ to $f(x)$, $f(x)$ to $x$ to $y_1$, $y_1$ to $z$, $z$ to $y_2$, and $y_2$ to $x$, see the right drawing of Figure
  \ref{fig_situation_weak_reduction_4}.

  Now we analyze the distance $d_H(a,b)$ between
  two vertices in $V(G)$. If $a$ and $b$ are both in $\tilde{C}_2$, then we can consider a shortest path $P$ in $H_2$. It may
  happen that $P$ uses some of the redirected edges. In this case $P$ contains at least two vertices from $\{x,y_1,y_2,z\}$.
  If $P$ uses more than two vertices from $\{x,y_1,y_2,z\}$ then we only consider those two vertices which have the largest distance
  on $P$. Looking at our redirected edges in $H$ we see, the distance between two such vertices is at most three, so that we have
  $d_H(a,b)\le d_{H_2}(a,b)+3$ in this case.

  Now let $b$ be in $\tilde{C}_2$. We consider a shortest path $P$ in $H_2$ from $z$ to
  $b$. In $H$ we have $d_H(f(x),z)\le 3$ by considering the path $[f(x),x,y_1,z]$. Since $d_H(z,y_2)=1$ we have
  $d_H(f(x),b)\le d_{H_2}(z,b)+4$. Similarly we obtain $d_H(w,b)\le d_{H_2}(z,b)+5$. With $D_2=D\backslash\{f(x)\}$ the
  set $D_2$ is a dominating set of $\tilde{C}_2$ and we can check that $|D_2|=\gamma(\tilde{C}_2)$ holds. Since $z\in D_2$ and
  $H_2$ is an optimal orientation, for
  $b_1\in D_2$, $b_2\notin D_2$ we have $d_{H_2}(z,b_1)\le \Xi(\gamma-1)-4$ and $d_{H_2}(z,b_2)\le \Xi(\gamma-1)-2$ yielding
  $d_H(f(x),b_1)\le \Xi(\gamma-1)$, $d_H(f(x),b_2)\le \Xi(\gamma-1)+2$, $d_H(w,b_1)\le \Xi(\gamma-1)+1$, and $d_H(w,b_2)\le
  \Xi(\gamma-1)+3$. This is compatible with $\Xi(\gamma)\le\Xi(\gamma-1)+4$ due to $f(x),b_1\in D$ and $w,b_2\notin D$.

  Now let $a$ be in $\tilde{C}_2$. we consider
  a shortest path $P$ in $H_2$ from $a$ to $z$. In $H$ we have $d_H(z,f(x))\le 4$ by considering the path $[z,y_2,x,w,f(x)]$. Since
  $P$ can not use an arc from $y_1$ to $z$ (this arc is directed in the opposite direction in $H_2$) either $P$ contains a vertex
  in $\{x,y_2\}$ or $P$ also exists in $H$, so that we have $d_H(a,f(x))\le d_{H_2}(a,z)+4$. Similarly we obtain $d_H(a,w)\le
  d_{H_2}(a,z)+3$. Since $H_2$ we conclude similarly as in the above paragraph that all distances are compatible with
   $\Xi(\gamma)\le\Xi(\gamma-1)+4$.
\end{Proof}

\begin{Lemma}
  \label{lemma_is_cut_vertex}
  Let $G$ be a minimal subgraph of a pair $(G', D)$ in first standard form. If there
  exist $z_1,z_2\in V(G)\backslash D$ with $f(z_1)=f(z_2)$ and $\{z_1,z_2\}\in E(G)$, then either $z_1$ or $z_2$ is a cut vertex.
\end{Lemma}
\begin{Proof}
  If $z_1$ has no other neighbors besides $z_2$ and $x:=f(z_1)$ then either $z_2$ is a cut vertex or $z_1$ can be deleted from $G$
  without destroying the properties of Definition \ref{def_minimal_subgraph}. We assume that whether $z_1$ nor $z_2$ is a cut
  vertex. Thus both $z_1$ and $z_2$ have further neighbors $y_1$ and $y_2$, respectively. Since $\{z_1,z_2\}$ can not be deleted
  we have $y_1\neq y_2$. Let $P_1$ be a shortest path from $y_1$ to $z_2$ in $G\backslash\{z_1\}$. Since $\{z_1,z_2\}$ can not be
  deleted $P_1$ contains the edge $\{x,z_2\}$. Similarly there exists a shortest path from $y_2$ to $z_1$ containing the edge
  $\{x,z_1\}$. Thus in the end the existence of $P_1$ and $P_2$ shows that $\{z_1,z_2\}$ could be deleted, which is a
  contradiction to the minimality of $G$.
\end{Proof}

\begin{Lemma}
  \label{lemma_no_quadrangle}
  Let $G$ be a minimal subgraph of a pair $(G', D)$ in first standard form. Let $x,y_1,y_2$ be three vertices not in the
  dominating set $D$ with $\{x,y_1\},\{x,y_2\}\in E(G)$ and $f(y_1)\neq f(x)\neq f(y_2)$ either one vertex of $x$, $y_1$,
  $y_2$ is a cut vertex, or $f(y_1)\neq f(y_2)$.
\end{Lemma}
\begin{Proof}
  We assume as contrary that none of $x$, $y_1$, $y_2$ is a cut vertex and $f(y_1)=f(y_2)$. Now we consider $G\backslash\{x\}$,
  which must be connected. Thus there must exist a path $P$ connecting $f(x)$ to $f(y_1)=f(y_2)$ and either one of the edges
  $\{x,y_1\}$, $\{x,y_2\}$ is a chord or one of the vertices $y_1$, $y_2$ could be deleted from $G$, which is a contradiction
  to the minimality of $G$.
\end{Proof}

\section{Proof of the main theorem}

In this section we want to prove Theorem \ref{thm_four_gamma}. We use the techniques of induction on $\gamma(G)$ and minimal
counter examples with respect to $\gamma(G)$.

\begin{Definition}
  We call a minimal subgraph $G$ of $(G',D)$ in first standard form a minimal counter example to Theorem \ref{thm_four_gamma}
  if we have $\max\Big\{diam_0(H,D)+4,diam_1(H,D)+2,diam_2(H,D)\Big\}>4\gamma$ for a minimal orientation $H$ and $\gamma=|D|$ is
  minimal with this property.
\end{Definition}

\begin{figure}[htp]
\begin{center}
  \setlength{\unitlength}{1cm}
  \begin{picture}(3.4,3.8)
    \put(0.2,0.2){\circle*{0.4}}
    \put(0.2,1.2){\circle{0.4}}
    \put(0.2,2.2){\circle{0.4}}
    \put(0.2,3.2){\circle*{0.4}}
    \put(1.2,0.2){\circle{0.4}}
    \put(2.2,0.2){\circle{0.4}}
    \put(3.2,0.2){\circle*{0.4}}
    \put(3.2,3.2){\circle*{0.4}}
    \put(3.2,2.2){\circle{0.4}}
    \put(3.2,1.2){\circle{0.4}}
    \put(2.2,3.2){\circle{0.4}}
    \put(1.2,3.2){\circle{0.4}}
    \put(0.4,0.2){\vector(1,0){0.6}}
    \put(1.4,0.2){\vector(1,0){0.6}}
    \put(2.4,0.2){\vector(1,0){0.6}}
    \put(3.2,0.4){\vector(0,1){0.6}}
    \put(3.2,1.4){\vector(0,1){0.6}}
    \put(3.2,2.4){\vector(0,1){0.6}}
    \put(3.0,3.2){\vector(-1,0){0.6}}
    \put(2.0,3.2){\vector(-1,0){0.6}}
    \put(1.0,3.2){\vector(-1,0){0.6}}
    \put(0.2,3.0){\vector(0,-1){0.6}}
    \put(0.2,2.0){\vector(0,-1){0.6}}
    \put(0.2,1.0){\vector(0,-1){0.6}}
    \put(0.05,-0.25){$v_0$}
    \put(1.05,-0.25){$v_1$}
    \put(2.05,-0.25){$v_2$}
    \put(3.05,-0.25){$v_3$}
    \put(3.45,1.10){$v_4$}
    \put(3.45,2.10){$v_5$}
    \put(3.05,3.55){$v_6$}
    \put(2.05,3.55){$v_7$}
    \put(1.05,3.55){$v_8$}
    \put(0.05,3.55){$v_9$}
    \put(-0.55,2.10){$v_{10}$}
    \put(-0.55,1.10){$v_{11}$}
  \end{picture}
  \quad\quad\quad\quad\quad\quad
  \begin{picture}(6.4,3.8)
    \put(0.2,0.2){\circle*{0.4}}
    \put(0.2,1.2){\circle{0.4}}
    \put(1.2,0.2){\circle{0.4}}
    \put(1.2,2.2){\circle{0.4}}
    \put(0.2,2.2){\circle*{0.4}}
    \put(0.2,3.2){\circle{0.4}}
    \put(2.2,2.2){\circle{0.4}}
    \put(2.2,3.2){\circle*{0.4}}
    \put(3.2,3.2){\circle{0.4}}
    \put(4.2,2.2){\circle{0.4}}
    \put(5.2,1.2){\circle{0.4}}
    \put(5.2,0.2){\circle{0.4}}
    \put(6.2,0.2){\circle*{0.4}}
    \put(6.2,1.2){\circle{0.4}}
    \put(5.2,2.2){\circle*{0.4}}
    \put(4.2,0.2){\circle{0.4}}
    \put(3.2,0.2){\circle*{0.4}}
    \put(2.2,0.2){\circle{0.4}}
    \put(1.2,0.2){\circle{0.4}}
    \put(1.4,0.2){\vector(1,0){0.6}}
    \put(2.4,0.2){\vector(1,0){0.6}}
    \put(3.4,0.2){\vector(1,0){0.6}}
    \put(4.4,0.2){\vector(1,0){0.6}}
    \put(5.2,0.4){\vector(0,1){0.6}}
    \put(5.2,1.4){\vector(0,1){0.6}}
    \put(5.0,2.2){\vector(-1,0){0.6}}
    \put(4.0,2.2){\vector(-1,0){1.6}}
    \put(2.0,2.2){\vector(-1,0){0.6}}
    \put(1.2,2.0){\vector(0,-1){1.6}}
    \put(1.0,0.2){\vector(-1,0){0.6}}
    \put(0.2,0.4){\vector(0,1){0.6}}
    \put(0.335,1.065){\vector(1,-1){0.730}}
    \put(1.0,2.2){\vector(-1,0){0.6}}
    \put(0.2,2.4){\vector(0,1){0.6}}
    \put(0.335,3.065){\vector(1,-1){0.730}}
    \put(2.2,2.4){\vector(0,1){0.6}}
    \put(2.4,3.2){\vector(1,0){0.6}}
    \put(3.065,3.065){\vector(-1,-1){0.730}}
    \put(5.4,0.2){\vector(1,0){0.6}}
    \put(6.2,0.4){\vector(0,1){0.6}}
    \put(6.065,1.065){\vector(-1,-1){0.730}}
    \put(1.05,-0.25){$v_8$}
    \put(2.05,-0.25){$v_9$}
    \put(3.05,-0.25){$v_0$}
    \put(4.05,-0.25){$v_1$}
    \put(5.05,-0.25){$v_2$}
    \put(4.60,1.10){$v_3$}
    \put(-0.50,1.10){$w_8$}
    \put(-0.40,0.10){$z_8$}
    \put(-0.50,3.10){$w_7$}
    \put(-0.40,2.10){$z_7$}
    \put(5.45,2.10){$v_4$}
    \put(6.45,0.10){$z_2$}
    \put(6.45,1.10){$w_2$}
    \put(1.05,2.55){$v_7$}
    \put(2.05,1.75){$v_6$}
    \put(4.05,2.55){$v_5$}
    \put(2.05,3.55){$z_6$}
    \put(3.05,3.55){$w_6$}
  \end{picture}
  \caption{The situation of Lemma \ref{lemma_no_circle} and the situation of Lemma \ref{lemma_no_special_circle}.}
  \label{fig_situation_no_circle}
\end{center}
\end{figure}
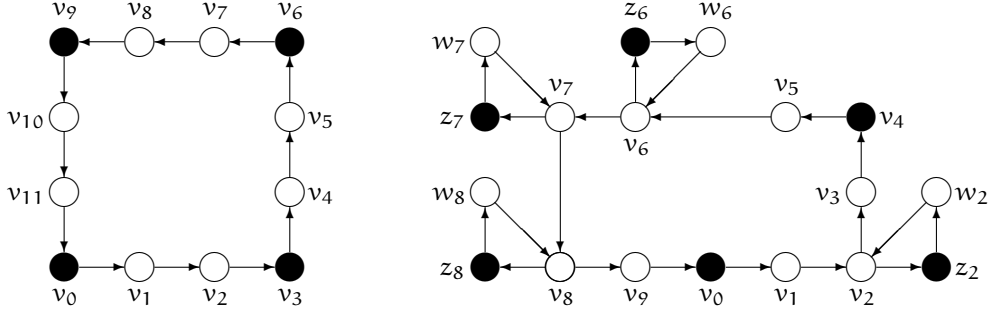

\begin{Lemma}
  \label{lemma_no_circle}
  Let $G$ be a minimal subgraph of $(G',D)$ in first standard form which is a minimal counter example to Theorem
  \ref{thm_four_gamma}, then there can not exist an elementary cycle  $C=[v_0,\dots,v_{3k}=v_0]$ in $G$ with $k\ge 2$ and
  the $v_{3j}\in D$ for all $0\le j< k$.
\end{Lemma}
\begin{Proof}
  We assume the existence of such a cycle $C$, see the left graph in Figure \ref{lemma_no_circle} for an example, and
  consider another graph $\tilde{G}$ arising from $G$ by:
  \begin{itemize}
    \item[(1)] deleting the edges of $C$,
    \item[(2)] deleting the vertices $v_{3j}$ for $0<j<k$,
    \item[(3)] inserting vertices $u_j$ and edges $\{v_0,v_j\}$, $\{v_0,u_j\}$, $\{u_j,v_j\}$ for all $0<j<3k$ with $3\nmid j$,
               and by
    \item[(4)] identifying all vertices $v_{3j}\in G$ with the vertex $v_0\in\tilde{G}$, meaning that we replace
               edges $\{v_{3j},x\}$ in $G$ by edges $\{v_0,x\}$ in $\tilde{G}$.
  \end{itemize}
  We remark that this construction does not produce multiple edges since $(G',D)$ is in first standard form.
  The set $\tilde{D}:=D\backslash\Big\{v_3,v_6,\dots,v_{3k-3}\Big\}$ is a dominating set of $\tilde{G}$ with
  $|\tilde{D}|=|D|-k+1$. Let
  $\tilde{H}$ be an minimal orientation of $(\tilde{G},\tilde{D})$. We construct an orientation $H$ of $G$ by taking over the
  directions of all common edges with $\tilde{H}$ and by orienting the edges of $C$ from $v_j$ to $v_{j+1}$, see the left graph
  in Figure \ref{lemma_no_circle}.

  Now we analyze the distances in $H$. For brevity we set $I:=\Big\{v_{3j}\,:\, 0\le j< k\Big\}$ (these are the vertices in $G$ which
  are associated with $v_0$ in $\tilde{G}$). The distance of two vertices in $I$ in the orientation $H$ is at most $3k-3$ and the
  distance of two vertices in $V(C)$ is at most $3k-1$. Thus we may assume $|D|> k$. Let $a,b$ be vertices in $V(G)$.
  \begin{itemize}
    \item[(1)] If $a$ and $b$ are elements of $\{v_j\,:\,0\le j< 3k\}$ then we have $d_H(a,b)\le 3k-1< 4|D|-4$.
    \item[(2)] If $a$ and $b$ are not in $I$ then we consider a shortest path $\tilde{P}$ in $\tilde{H}$ connecting
               $a$ and $b$.
    \item[(3)] If $a\in I$ and $b\notin I$ then we consider a shortest path $\tilde{P}$ in $\tilde{H}$ connecting
               $v_0$ and $b$.
    \item[(4)] The case $a\notin I$ and $b\in I$  then we consider a shortest path $\tilde{P}$ in $\tilde{H}$ connecting
               $a$ and $v_0$.
  \end{itemize}
  Let $\tilde{P}$ be an arbitrary shortest path in $\tilde{H}$ connecting $a$ and $b$. It may happen that in $H$ this path
  $\tilde{P}$ does not exist since it may contain the vertex $v_0$ corresponding to two different vertices $v_{3i}$ and $v_{3j}$
  in $G$ or may contain one of the edges $\{v_0,v_j\}$, $\{v_0,u_j\}$, or $\{u_j,v_j\}$ with $3\nmid j$.

  Now we want to construct a path $P$ which does connect $a$ and $b$ in $H$. The path $\tilde{P}$ may use one of the edges
  $\{v_0,v_j\}$, $\{v_0,u_j\}$, or $\{u_j,v_j\}$ with $3\nmid j$. Deleting all these edges decomposes $\tilde{P}$ in at least
  two parts $\tilde{P}_1$, $\dots$, $\tilde{P}_m$ with $|\tilde{P}_1|+|\tilde{P}_m|\le |\tilde{P}|-1$. Using a suitable segment
  $\tilde{C}$ of the cycle $C$ we obtain a path $P=\tilde{P}_1\cup\tilde{C}\cup\tilde{P}_m$ of length at most
  $|\tilde{P}_1|+|\tilde{P}_m|+|\tilde{C}| \le |\tilde{P}|+3k-2$. If $\tilde{P}$ does not use one of these edges then it can
  only happen that $v_0$ is used in $\tilde{P}$ corresponding to two different vertices $v_{3i}$ and $v_{3j}$ in $G$.
  In this case we can use a suitable segment $\tilde{C}$ of the cycle $C$, which starts and ends in a vertex of $I$, to obtain
  a path $P$ connecting $a$ and $b$ in $H$ of length at most $|\tilde{P}|+3k-3$.

  Now we are ready to prove that $G$ is not a counter example. If $\gamma(\tilde{G})<|\tilde{D}|$ then we have
  $diam\left(\tilde{H}\right)
  \le 4\cdot|\tilde{D}|-4=4\cdot |D|-4k$ due to the minimality of $G$. In each of the cases (1)-(4)
  we have $d_H(a,b)\le 4\cdot |D|-k-2\le 4\cdot |D|-4$ for all $a,b\in G$. Otherwise we have $\gamma(\tilde{G})=|\tilde{D}|$ and
  $\tilde{D}$ is a minimal dominating set of $\tilde{G}$. In this case we have
  \begin{eqnarray*}
    diam_2(H,D) &\le & \max\Big\{diam_2\left(\tilde{H},\tilde{D}\right)+3k-2,diam_1\left(\tilde{H},\tilde{D}\right)+3k-1,3k-1\Big\}\\
                &\le & 4\cdot |D|-k+2\\
                &\le & 4\cdot|D|\\
    diam_1(H,D) &\le & \max\Big\{diam_1\left(\tilde{H},\tilde{D}\right)+3k-2,diam_0\left(\tilde{H},\tilde{D}\right)+3k-1,3k-1\Big\}\\
                &\le & 4\cdot |D|-k\\
                &\le & 4\cdot|D|-2\\
    diam_0(H,D) &\le & \max\Big\{diam_0\left(\tilde{H},\tilde{D}\right)+3k-2,3k-3\Big\}\\
                &\le & 4\cdot |D|-k-2\\
                &\le & 4\cdot|D|-4\\
  \end{eqnarray*}
\end{Proof}

\begin{Lemma}
  \label{lemma_no_special_circle}
  Let $G$ be a minimal subgraph of $(G',D)$ in first standard form which is a minimal counter example to Theorem
  \ref{thm_four_gamma}, then there can not exist an elementary cycle  $C=[v_0,\dots,v_{l}=v_0]$ in $G$ with the
  following properties:
  \begin{itemize}
   \item[(1)] $v_0\in D$,
   \item[(2)] $|V(C)\cap D|\ge 2$,
   \item[(3)] $l\ge 6$, and
   \item[(4)] if $v_j\notin D$ then either $f(v_j)\in \{v_{j-1},v_{j+1}\}$ or $v_j$ is a cut vertex in $G$ where the component
              containing $f(v_j)$ contains exactly one vertex of $D$.
  \end{itemize}
\end{Lemma}
\begin{Proof}
  We assume the existence of such a cycle $C$. By $y$ we denote the number of cut vertices $v_j$ in $C$ and by $Y$
  the corresponding set. For all $v\in Y$ we have $f(v)\notin C$ since otherwise we could apply Lemma \ref{lemma_weak_reduction_3}.
  If $e=\{v',v''\}$ would be a chord of $C$ then $|\{v',v''\}\cap D|=1$ since $(G',D)$ is in first standard form and $G$ is
  a minimal subgraph, which especially means that we can not delete the edge $e$. We assume w.l.o.g. $v'\in D$ and conclude
  $f(v'')=v'$. Thus $v''$ is not a cut vertex and due to property (4) the edge $e$ is not a chord. Finally we conclude that $C$
  is chordless. For $y=0$ we would have $v_{3j}\in D$ due to $l\ge 6$ and the property $f(v_j)\in \{v_{j-1},v_{j+1}\}$ for
  vertices $v_j\notin D$. Thus we may assume $y\ge 1$ since otherwise we could apply Lemma \ref{lemma_no_circle}. For each
  $v_j\in Y$ we set $z_j=f(v_j)\notin V(C)$ and denote by $w_j\in V(G)\backslash(V(C)\cup D)$ the vertex which is adjacent
  to $v_j$ and $z_j$. By $k$ we denote the number of vertices $v_j$ in $V(C)$
  which are also contained in $D$. 
  Due to condition (2) we have $k\ge 2$. The two neighbors on the cycle $C$ of a vertex in $Y$ both are not contained in $D$.
  For a vertex $v\in V(C)\backslash(D\cup Y)$ one neighbor on $C$ is $f(v)$ and the other neighbor lies in $V(C)\backslash D$.
  Thus the length $|C|$ of the cycle is given by $3k+y\ge 7$. On the right hand side of Figure \ref{lemma_no_circle} we have
  depicted an example with $k=2$ and $y=4$.

  Now we consider another graph $\tilde{G}$ arising from $G$ by:
  \begin{itemize}
    \item[(1)] deleting the edges of $C$,
    \item[(2)] deleting the vertices $\Big(\left\{z_j,w_j\,:\,0<j<l\right\}\cup(V(C)\cap D)\Big)\backslash\left\{v_0\right\}$,
    \item[(3)] inserting vertices $u_j$ and edges $\{v_0,v_j\}$, $\{v_0,u_j\}$ ,$\{u_j,v_j\}$ for all $0<j< l$ with $v_j\notin D$,
               and by
    \item[(4)] identifying all vertices $v_{j}\in D$ with the vertex $v_0\in\tilde{G}$, meaning that we replace
               edges $\{v_j,x\}$ in $G$ by edges $\{v_0,x\}$ in $\tilde{G}$.
  \end{itemize}
  We remark that this construction does not produce multiple edges since $(G',D)$ is in first standard form.
  The set $\tilde{D}:=D\backslash\Big\{v_1\dots,v_{l-1},z_1,\dots,z_{l-1}\Big\}$ is a dominating set of $\tilde{G}$ with
  $|\tilde{D}|=|D|-k-y+1$. Let $\tilde{H}$ be an minimal orientation of $(\tilde{G},\tilde{D})$. We construct an
  orientation $H$ of $G$ by taking over the directions of all common edges with $\tilde{H}$ and by orienting the edges
  of $C$ from $v_j$ to $v_{j+1}$. The missing edges corresponding to $z_j$ and $w_j$ are oriented from $v_j$ to $z_j$,
  from $z_j$ to $w_j$, and from $w_j$ to $v_j$, see the graph on the right hand side of Figure \ref{lemma_no_circle}.
  For brevity we set $A=V(C)\cup\Big\{w_j,z_j\,:\, 0<j<l\Big\}$.

  Now we analyze the distances in $H$. For $a_1,b_1\in A$ we have $d_H(a_1,b_1)\le 3k+y+3$, for $a_2,b_2\in V(C)$ we
  have $d_H(a_2,b_2)\le 3k+y-1$, and for $a_3,b_3\in V(C)\cap D$ we have
  $d_H(a_3,b_3)\le 3k+y-3$. Thus we may assume $|D|> k+y$. Let $a,b$ be vertices in $V(G)$.
  \begin{itemize}
    \item[(1)] If $a$ and $b$ are elements of $A$ then we have $d_H(a,b)\le 3k+y+3< 4|D|-4$.
    \item[(2)] If $a$ and $b$ are not in $A$ then we consider a shortest path $\tilde{P}$ in $\tilde{H}$ connecting
               $a$ and $b$.
    \item[(3)] If $a\in A$ and $b\notin A$ then we consider a shortest path $\tilde{P}$ in $\tilde{H}$ connecting
               $v_0$ and $b$.
    \item[(4)] The case $a\notin A$ and $b\in A$  then we consider a shortest path $\tilde{P}$ in $\tilde{H}$ connecting
               $a$ and $v_0$.
  \end{itemize}
  Let $\tilde{P}$ be a shortest path in $\tilde{H}$ connecting two vertices $a$ and $b$. Similarly as in the proof
  of Lemma \ref{lemma_no_circle} we construct a path $P$ in $H$ connecting $a$ and $b$. Doing the same analysis we obtain
  $|P|\le |\tilde{P}|+3k+y-2$. Starting or ending at a vertex $z_i$ or $w_i$ increases the length by at most $2$.

  If $\gamma\left(\tilde{G}\right)<|\tilde{D}|=|D|-k-y+1$ then we would have $d_H(a,b)\le 4|D|-k-3y+1\le 4|D|-4$. Thus
  we may assume $\gamma\left(\tilde{G}\right)=|\tilde{D}|= |D|-k-y+1$,
  meaning that $\tilde{D}$ is a minimal dominating set. With this clearly we have $d_{\tilde{H}}(v_0,b),d_{\tilde{H}}(a,v_0)\le
  4\cdot |\tilde{D}|-2$ for all $a,b\in \tilde{G}$ and $d_{\tilde{H}}(v_0,b'),d_{\tilde{H}}(a',v_0)\le 4\cdot |\tilde{D}|-4$ for
  all $a',b'\in\tilde{D}$.

  For $k+y\ge 3$ and $|D|\ge k+y+1$ we have
  \begin{eqnarray*}
    diam_2(H,D) &\le & \max\Big\{diam_2\left(\tilde{H},\tilde{D}\right)+3k+y-2,diam_1\left(\tilde{H},\tilde{D}\right)
                 +3k+y,3k+y+3\Big\}\\
                &\le & 4\cdot |D|-k-3y+2\\
                &\le & 4\cdot|D|\\
    diam_1(H,D) &\le & \max\Big\{diam_1\left(\tilde{H},\tilde{D}\right)+3k+y,diam_0\left(\tilde{H},\tilde{D}\right)
                 +3k+y,3k+y+3\Big\}\\
                &\le & 4\cdot |D|-k-3y+2\\
                &\le & 4\cdot|D|-2\\
    diam_0(H,D) &\le & \max\Big\{diam_0\left(\tilde{H},\tilde{D}\right)+3k+y,3k+y+3\Big\}\\
                &\le & 4\cdot |D|-k-3y\\
                &\le & 4\cdot|D|-4\\
  \end{eqnarray*}
\end{Proof}

\noindent
Now we are ready to prove Theorem \ref{thm_four_gamma}:
\begin{Proof}(of Theorem \ref{thm_four_gamma})\\
  Let $G$ be a minimal subgraph of $(G',D)$ in first standard form which is a minimal counter example to Theorem
  \ref{thm_four_gamma}. Due to Lemma \ref{lemma_small_exact_values} and Lemma \ref{lemma_gamma_3} we can assume $|D|\ge 4$.
  We show that we have
  $|V(G)|\le 4\cdot(|D|-1)+1$. In this case we can utilize an arbitrary
  orientation $H$ of $G$. Since a shortest path uses every vertex at most once we would have $diam(H)\le 4\cdot(|D|-1)$.
  Applying Lemma \ref{lemma_subgraph} we conclude $diam_{min}\left(G'\right)\le 4\cdot|D|=4\cdot\gamma(G')$, which is a contradiction
  to $G$ being a minimal counter example to Theorem \ref{thm_four_gamma} and instead proves this theorem.

  At first we summarize some structure results for minimal counter examples to Theorem \ref{thm_four_gamma}.
  \begin{itemize}
    \item[(1)] We can not apply one of the Lemmas \ref{lemma_weak_reduction_2}, \ref{lemma_reduction_3}, \ref{lemma_reduction_5},
               \ref{lemma_weak_reduction_3}, or \ref{lemma_is_cut_vertex}. So if $v\in V(G)$ is a cut vertex we have $v\notin D$
               and there exists a unique vertex $t(V)\notin D$ such that we have $\{v,f(v)\},\{f(v),t(v)\},\{t(v),v\}\in E(G)$ and
               all neighbors of $f(v),t(v)$ are contained in $\{f(v),t(v),v\}$.
    \item[(2)] Due to Lemma \ref{lemma_weak_reduction_4}, Lemma \ref{lemma_no_quadrangle}, and (1) there do not exist
               pairwise different vertices $x,y_1,y_2\in V(G)\backslash D$ with $\{x,y_1\},\{x,y_2\}\in E(G)$ and $f(y_1)=f(y_2)$.
    \item[(3)] We can not apply Lemma \ref{lemma_no_circle} or Lemma \ref{lemma_no_special_circle} on $G$.
  \end{itemize}
  In order to bound $|V(G)|$ from above we perform a technical trick and count the number of vertices of a different graph
  $\tilde{G}$. Therefore we label the cut vertices of $G$ by $v_1,\dots,v_m$. With this we set
  $$
    \tilde{D}=\Big(D\cup\left\{v_i\,:\,1\le i\le m\right\}\Big)\backslash\left\{f(v_i)\,:\,1\le i\le m\right\}.
  $$
  The graph $\tilde{G}$ arises from $G$ by deleting the $f(v_i), t(v_i)$ for $1\le i\le m$ and by replacing the remaining
  edges $\{v_i,x\}$ by a pair of two edges $\{v_i,y_{x,i}\},\{y_{x,i}\},x\}$, where the $y_{x,i}$ are new vertices. We
  have $|\tilde{D}|=|D|$, $|V(\tilde{G})|\ge|V(G)|$, the set $\tilde{D}$ is a dominating set of $\tilde{G}$, and $\tilde{G}$
  is a subgraph of a suitable pair in first standard form. If $\tilde{G}$ would not be a minimal subgraph than also $G$ would
  not be a minimal subgraph. We have the following structure results for $\tilde{G}$:
  \begin{itemize}
    \item[(a)] There do not exist two vertices $u,v\in V(\tilde{G})\backslash\tilde{D}$ with $\{u,v\}\in E(\tilde{G})$ and
               $f(u)=f(v)$.
    \item[(b)] There do not exist pairwise different vertices $x,y_1,y_2\in
               V(\tilde{G})\backslash \tilde{D}$ with $\{x,y_1\},\{x,y_2\}\in E(\tilde{G})$ and $f(y_1)=f(y_2)$.
    \item[(c)] We can not apply Lemma \ref{lemma_no_circle} or Lemma \ref{lemma_no_special_circle} on $\tilde{G}$.
  \end{itemize}
  Since our construction of $\tilde{G}$ has removed all such configurations (a) holds. If in (b) $f(y_1)=f(y_2)$ is an
  element of $D$ then such a configuration also exists in $G$, which is a contradiction to (2). If $f(y_1)$ corresponds
  to a $v_i$ in $G$, then $y_1$ and $y_2$ would correspond to two new vertices $y_{i,e}$ and $y_{i,e'}$. In this case
  we would have a double edge from $x$ to $v_i$ in $G$, which is not true. Thus (b) holds. Since all vertices in
  $\tilde{D}\backslash D$ correspond to cut vertices in $G$ also (c) holds.

  In order to prove $|V(\tilde{G})|\le 4\cdot(|\tilde{D}|-1)+1$ we construct a tree $T$ fulfilling
  \begin{itemize}
    \item[(i)] $\tilde{D}\subseteq V(T)$ and
    \item[(ii)] if $v_1\in V(T)\backslash \tilde{D}$ then we have $\{f(v_1),v_1\}\in E(T)$.
  \end{itemize}
  Therefore we iteratively construct trees $T_k$ for $1\le k\le |\tilde{D}|$. The tree $T_1$ is composed of a single vertex
  $x_1\in \tilde{D}$. The tree $T_1$ clearly fulfills condition (ii). To construct $T_{k+1}$ from $T_k$ we find a
  vertex $x_{k+1}$ in $\tilde{D}\backslash V(T_k)$ with the minimum distance to $T_k$. The tree $T_{k+1}$ is the union of
  $T_k$ with a shortest path $P_{k+1}$ from $x_{k+1}$ to $T_k$. Since $\tilde{D}$ is a dominating set this path $P_{k+1}$
  has length at most three. Since $\tilde{G}$ is a subgraph of a suitable pair in first standard form $P_{k+1}$ has length
  at least two. For $P_{k+1}=[x_{k+1},v_1,v_2]$ we have $v_1,v_2\notin \tilde{D}$ due to the first standard form and
  $f(v_1)=x_{k+1}$, $v_2\in V(T_k)$. Since condition (ii) is fulfilled for $T_k$ it is also fulfilled for $T_{k+1}$ in this
  case. In the remaining case we have
  $P_{k+1}=[x_{k+1},v_1,v_2,v_3]$ with $v_1,v_2\notin V(T_k)$, $v_1,v_2\notin \tilde{D}$, and $v_3\in V(T_k)$.
  If $f(v_2)$ would not be contained in $V(T_k)$ then $[f(v_2),v_2,v_3]$ would be a shorter path connecting $f(v_2)$ to $T_k$.
  Thus we have $f(v_2)\in V(T_k)$ and we may assume $v_3=f(v_2)$. (We may simply consider the path $[x_{k+1},v_1,v_2,f(v_2)]$
  instead of $P_{k+1}$.) Due to $x_{k+1}\notin V(T_k)$ and $T_k$ fulfilling condition (ii), these conditions are also
  fulfilled for $T_{k+1}$. In the end we obtain a tree $T_{|\tilde{D}|}$ fulfilling condition (i) and condition (ii). By
  considering the paths $P_k$ we conclude $|V(T)|\le |\tilde{D}|+2(|\tilde{D}|-1)$.

  \medskip

  Clearly we have some alternatives during the construction of $T_{|\tilde{D}|}$. Now we assume that $T$ is a subtree of
  $\tilde{G}$ fulfilling conditions (i) and (ii), and having the maximal number of vertices. In the next step we want to
  prove some properties of the vertices in $T$.

  Let $v\in \tilde{D}$ and let $u\in V(\tilde{G})\backslash V(T)$ be a neighbor of $v$ in $\tilde{G}$. We prove that every neighbor
  $u'$ of $u$ in $\tilde{G}$ is contained in $V(T)$. Clearly we have $u'\notin \tilde{D}$. Due to (a) we have $f(u')\neq v$.
  If  $u'\notin V(T)$ then adding the edges
  $A:=\Big\{\{v,u\},\{u,u'\},\{u',f(u')\}\Big\}$ gives an elementary cycle $C=[v_0,\dots,v_l]$ in
  $\Big(V(T)\cup\{u,u'\},E(T)\cup A\Big)$, where $v_0=v_l$ and $l\ge 6$. Since we can not apply Lemma \ref{lemma_no_circle}
  there exists an index $j$ (reading the indices modulo $l$) fulfilling
  $$
    v_j\in \tilde{D}\quad\text{and}\quad v_{j+1},v_{j+2},v_{j+3}\in V(T)\backslash \tilde{D}.
  $$
  Since the edge $\{v_{j+1},v_{j+2}\}$ is contained in $E(T)$ also the edge $\left\{v_{j+2},f\left(v_{j+2}\right)\right\}$
  is contained in $E(T)$. Similarly we conclude that the edge $\{v_{j+3},f(v_{j+3})\}$ is contained in $E(T)$. If $v_{j+1}$
  has no further neighbors besides $v_j$ and $v_{j+2}$ in $T$ then
  $$
    T':=\left((V(T)\cup\{u,u'\})\backslash\{v_{j+1}\},(E(T)\cup A)\backslash
    \Big\{\{v_j,v_{j+1}\},\{v_{j+1},v_{j+2}\}\Big\}\right)
  $$
  would be a subtree of $\tilde{G}$ fulfilling the conditions (i) and (ii) with a larger number of vertices than $T$.
  Thus such an $u'$ can not exist in this case. If $v_{j+1}$ has further neighbors in $T$, then deleting the edge
  $\{v_{j+1},v_{j+2}\}$ and adding the edges and vertices of $A$ would also yield a subtree of $\tilde{G}$ fulfilling the
  conditions (i) and (ii) with a larger number of vertices than $T$.

  The same statement also holds for $v\in V(T)\backslash\tilde{D}$ since we may consider $f(u)$ instead $v$. Thus in $\tilde{G}$
  we have $\{u,v\}\cap V(T)\neq\emptyset$ for every edge $\{u,v\}\in E(\tilde{G})$.

  \medskip

  For a graph $K$ and a vertex $v\in V(K)$ we denote by $S(K,v)$ the uniquely defined maximal bridgeless connected subgraph of $K$
  containing $v$. If every edge being adjacent to $v$ is a bridge or $v$ do not have any edges, then $S$ consists only of vertex $v$.
  We remark that $u\in S(K,v)$ is an equivalence relation ${\sim}{_K}$ for all vertices $u,v\in V(K)$.
  By $F$ we denote the set of vertices in $V(T)$ which are either contained in $\tilde{D}$ or have a degree in $V(T)$ of at least
  three. We have
  $$
    |V(T)|+|F|\le 4\cdot|\tilde{D}|-2,
  $$
  which can be proved by induction on $|V\left(T_k\right)|+|F\cap V\left(T_k\right)|\le 4\cdot k-2$ for $1\le k\le |\tilde{D}|$.
  Clearly we have $|V\left(T_1\right)|+|F\cap V\left(T_1\right)|=2\le 4\cdot 2-2$. The tree $T_{k+1}$ arises from $T_k$ by adding
  a path $P_{k+1}$ of length at most three. If $|P_{k+1}|=3$ then we have $F\cap V\left(T_{k+1}\right)=\Big(F\cap
  V\left(T_{k}\right)\Big)\cup\{x_{k+1}\}$ and $|V\left(T_{k+1}\right)|\le|V\left(T_k\right)|+3$. For $|P_{k+1}|=2$ we have
  $|V\left(T_{k+1}\right)|\le|V\left(T_k\right)|+2$ and $|F\cap V\left(T_{k+1}\right)|\le|F\cap V\left(T_k\right)|+2$.

  For a graph $K$ containing $T$ as a subgraph we denote by $N(K)$ the number $\Big|\{S(K,v)\,:\,v\in F\}\Big|$ of
  equivalence classes of $\sim_K$. Since $T$ is a tree we have $N(T)=|F|$. Now we recursively construct a sequence of graphs
  $G_i$ for $1\le i\le|F|$ fulfilling
  \begin{equation}
    \label{eq_condition}
    |V(G_i)|+N(G_i)\le 4\cdot|\tilde{D}|-2,\, N(G_i)\le i,\,\text{and}\, T\subseteq G_i\subseteq\tilde{G}.
  \end{equation}
  This yields a graph $G_1$ containing at most $4\cdot|\tilde{D}|-3$ vertices, where each two elements of $\tilde{D}$ are connected
  by at least two edge disjoint paths. So either we have $|V(\tilde{G})|\le 4\cdot|\tilde{D}|-3$ or $\tilde{G}$ and $G$ are not
  minimal subgraphs.

  During the following analysis we often delete a vertex $v$ or an edge $e$ from the tree $T$ in such a way that it decomposes
  in exactly two subtrees $T^1$ and $T^2$. Since $T$ contains no cut vertices there exists a path $M$ in $\tilde{G}$ without $v$ or
  without $e$ connecting $T^1$ and $T^2$. Since there does not exist an edge $\{u_1,u_2\}\in E\left(\tilde{G}\right)$ with
  $\{u_1,u_2\}\cap V(T)=\emptyset$ we have $|M|\le 2$ if $M$ is a shortest path. 

  For $G_{|F|}=T$ condition \ref{eq_condition} holds. Now for $i\ge 2$ let $G_i$ be given. If there exists a vertex
  $u\in V(\tilde{G})\backslash V(G_i)$ having neighbors $x,y\in V(G_i)$ with $S(G_i,x)\neq S(G_i,y)$ we define $G_{i-1}$
  by adding vertex $u$ and adding all edges, being adjacent with $u$ in $\tilde{G}$, to $G_i$. With this we have
  $|V(G_{i-1})|=|V(G_i)|+1$ and $N(G_{i-1})=N(G_i)-1$, so that condition \ref{eq_condition} is fulfilled for $G_{i-1}$.

  Now we deal with the cases where $i\ge 2$ and where such vertices $u$, $x$, $y$ do not exist. We use the setwise defined distance
  $$
    d_K(A,B):=\min\Big\{d_K(a,b)\,:\,a\in A,\,b\in B\Big\}.
  $$
  Now we choose $f_1,f_2\in F$ with $S(G_i,f_1)\neq S(G_i,f_2)$, where $d_{G_i}(S(G_i,f_1),S(G_i,f_2))$ is minimal.
  Clearly we have $1\le d_{G_i}(S(G_i,f_1),S(G_i,f_2))\le 3$. By $P_{f_1,f_2}$ we denote the corresponding shortest
  path connecting $S(G_i,f_1)$ with $S(G_i,f_2)$.

  If $|P_{f_1,f_2}|=[v_0,v_1]$ and the edge $\{v_0,v_1\}$ is not contained in $E(T)$, then we simply add this edge to $G_i$ to
  obtain $G_{i-1}$. So we may assume that $\{v_0,v_1\}\in E(T)$. Deleting $\{v_0,v_1\}$ in $T$ decomposes $T$ into two subtrees
  $T^1$ and $T^2$, where we assume w.lo.g. that $f_1\in V(T^1)$ and $f_2\in V(T^2)$. Due to
  $d_{\tilde{G}\backslash\{v_0,v_1\}}(T^1,T^2)\le 2$ we can obtain a graph $G_{i-1}$ adding add most one vertex, where
  $S(G_{i-1},f_1)=S(G_{i-1},f_2)$ holds.

  If $|P_{f_1,f_2}|=[v_0,v_1,v_2]$ and $v_1\notin V(T)$ the we can add $v_1$ and add all its edges to $G_i$ to obtain
  $G_{i-1}$. So we may assume $v_1\in V(T)$. If $\{v_0,v_1\}$ or $\{v_1,v_2\}$ would not be contained in $E(T)$, then we may
  simply add it to $G_i$, without increasing the number of vertices, and are in a case $|P_{f_1,f_2}|=1$. So we may assume
  $\{v_0,v_1\},\{v_1,v_2\}\in E(T)$. Due to $S(G_i,v_0)\neq S(G_i,v_1)\neq S(G_i,v_2)$ and the minimality of $P_{f_1,f_2}$
  we have $v_1\notin F$. Thus $v_1$ has degree two in $T$ and removing $v_1$ decomposes $T$ into two subtrees
  $T^1$ and $T^2$, where we assume w.lo.g. that $f_1\in V(T^1)$ and $f_2\in V(T^2)$. Since there does not exist a
  cut vertex in $\tilde{G}$ we have $d_{\tilde{G}\backslash\{v_1\}}(T^1,T^2)\le 2$ and we can obtain a graph $G_{i-1}$
  adding add most one vertex, where $S(G_{i-1},f_1)=S(G_{i-1},f_2)$ holds.

  The remaining case is $|P_{f_1,f_2}|=[v_0,v_1,v_2,v_3]$. Due to the minimality of $P_{f_1,f_2}$ we have $f(v_2)\in V(S(G_i,f_2))$
  and $f(v_1)\in V(S(G_i,f_1))$. Thus we may assume $v_0,v_3\in\tilde{D}$. Additionally we have $\{v_1,v_2\}\cap V(T)\neq\emptyset$.
  If $v_j\notin V(T)$ we may simply add $v_j$ and its edges to $G_i$ to obtain $G_{i-1}$. So we may assume $v_1,v_2\in V(T)$.
  W.l.o.g. we assume $\{v_1,v_2\}\in E(T)$. Otherwise there exists an edge $\{v_1,v_4\}\in E(T)$ with $v_4\neq v_0$ and
  we could choose $f_1=f(v_1)$, $f_2=f(v_4)$. The vertices $v_1$ and $v_2$ both have degree two in $T$. Deleting $v_1$ in
  $T$ gives two subtrees $T^1$ and $T^2$, where we can assume $v_0\in V(T^1)$ and $v_2\in V(T^2)$. Since there does not exist a
  cut vertex in $\tilde{G}$ we have $d_{\tilde{G}\backslash\{v_1\}}(T^1,T^2)\le 2$ and denote the corresponding shortest path by
  $R_1$. If $R_1=[r_0,r_1,r_2]$ does not end in $v_2$ then we could obtain $G_{i-1}$ by adding vertex $r_1$ and its edges to $G_i$.
  Similarly we may delete vertex $v_2$ to obtain a shortest path $R_2$ which ends in $v_1$. But in this case the edge $\{v_1,v_2\}$
  could be deleted from $\tilde{G}$, which is a contradiction to the minimality of $\tilde{G}$.
\end{Proof}

We remark that we conjecture that if $G$ is a critical minimal subgraph of a pair $(G',D)$ in first standard form then we always can apply one of the lemmas \ref{lemma_reduction_1}, \ref{lemma_reduction_2}, \ref{lemma_reduction_3}, \ref{lemma_reduction_5}, \ref{lemma_weak_reduction_1}, \ref{lemma_weak_reduction_2}, \ref{lemma_weak_reduction_3}, \ref{lemma_weak_reduction_4}, \ref{lemma_is_cut_vertex}, \ref{lemma_no_quadrangle}, \ref{lemma_no_circle}, or \ref{lemma_no_special_circle}.

We would like to remark that our reduction technique is constructive in the following sense: If we have a graph $G$ and a dominating set $D$, not necessarily a minimal dominating set of $G$, then we can construct an orientation $H$ of $G$ in polynomial time fulfilling $diam(H)\le 4\cdot |D|$: At first we apply the transformations of the poof of Lemma \ref{lemma_first_standard} to obtain a graph $\tilde{G}$, which fulfills conditions (1), (3)-(6) of Definition \ref{Def_first_standard} and where $D$ remains a dominating set. In the following we will demonstrate how to obtain an orientation $\tilde{H}$ of $\tilde{G}$ fulfilling $diam\left(\tilde{H}\right)\le 4\cdot |D|$. From such an orientation we can clearly reconstruct an orientation $H$ of $G$. Since Lemma \ref{lemma_subgraph} does not use the minimality of the dominating set $D$ we can restrict our consideration on a minimal subgraph $\hat{G}$ of $\tilde{G}$. Since none of the lemmas in Section \ref{sec_reductions} uses the
  minimality of the domination set $D$, we can apply all these reduction steps on $\hat{G}$. These steps can easily be reversed afterwards. The proofs of Lemma \ref{lemma_no_circle} and Lemma \ref{lemma_no_special_circle} have to be altered very slightly to guarantee a suitable reduction also in the case where $D$ is not minimal. (Here only the analysis is affected, not the construction.) We end up with a graph $\dot{G}$ with dominating set $\dot{D}$ (here $\dot{D}$ arises from $D$ by applying the necessary reduction steps). Since in the proof of Theorem \ref{thm_four_gamma} we show $\left|V\left(\dot{G}\right)\right|\le 4\cdot \left|\dot{D}\right|-3$ we can choose an arbitrary strong orientation and reverse all previous steps to obtain an orientation $H$ of $G$ with $diam(H)\le 4\cdot |D|$. We remark that all steps can be performed in polynomial time.

\section{Conclusion and outlook}
\noindent
In this article we have proven
$$
    \overset{\longrightarrow}{diam}_{min}(G)\le 4\cdot\gamma(G)
$$
for all bridgeless connected graphs and conjecture
$$
  \overset{\longrightarrow}{diam}_{min}(G)\le\left\lceil\frac{7\gamma(G)+1}{2}\right\rceil
$$
to be the true upper bound. Lemma \ref{lemma_gamma_3} shows that Theorem \ref{thm_four_gamma} is not tight for $\gamma=3$. Some of our reduction steps in Section \ref{sec_reductions} can also be used for a proof of Conjecture \ref{main_conj}. Key ingredients might be the lemmas \ref{lemma_no_circle} and \ref{lemma_no_special_circle}, which can be utilized as reductions for Conjecture \ref{main_conj} if $k+y$ is \text{large enough}. Figure \ref{fig_gamma_3} indicates several cases which can not be reduced so far.

Besides a proof of Conjecture \ref{main_conj} one might consider special subclasses of general graphs to obtain stronger bounds on the minimum oriented diameter. E.~g. for $C_3$-free graphs and $C_4$-free graphs we conjecture that the minimum oriented diameter is at most $3\cdot\gamma +c$.

\nocite{1046.05025,1055.05042,0603.05040,1064.05050,1050.05114,1025.05020,1022.68599}

\bibliographystyle{amsplain}
\bibliography{oriented_diameter}
\end{document}